%% file: generating_sets.tex
\documentclass{elsart}

\usepackage{amsmath}
\usepackage{amssymb}
\usepackage{amstext}
\usepackage{epsfig}
\usepackage{epic}
\usepackage{natbib}

\newcommand{\Graph}{\mathcal{G}}
\newcommand{\Lattice}{\mathcal{L}}
\newcommand{\Feasible}{\mathcal{F}}
\newcommand{\FourTiTwo}{\texttt{4ti2}}
\newcommand{\CP}{\mathcal{CP}}
\newcommand{\NF}{\mathcal{NF}}

\newcommand{\boproof}{\textbf{Proof.} }
\newcommand{\eoproof}{\qed}

\newcommand{\sigmabar}{{\bar{\sigma}}}
\newcommand{\taubar}{{\bar{\tau}}}

\DeclareMathOperator{\supp}{supp}
\DeclareMathOperator{\rank}{rank}

\begin{document}
\begin{frontmatter}
\title{Computing generating sets of lattice ideals\thanksref{support}}
\author[Magdeburg]{Raymond Hemmecke},
\ead{hemmecke@math.uni-magdeburg.de}
\author[CORE]{Peter N. Malkin}
\ead{malkin@core.ucl.ac.be}

\thanks[support]{This work was supported by the European TMR
network ADONET 504438.}
\address[Magdeburg]{Otto-von-Guericke-University Magdeburg, Germany}
\address[CORE]{CORE and INMA, Universit\'e catholique de Louvain, Belgium}

\begin{abstract}
In this article, we present a new algorithm for computing generating sets and
Gr\"obner bases of lattice ideals. In contrast to other existing methods, our
algorithm starts computing in projected subspaces and then iteratively lifts the
results back into higher dimensions, by using a completion procedure, until the
original dimension is reached. We give a completely geometric presentation of
our Project-and-Lift algorithm and describe also the two other existing main
algorithms in this geometric framework. We then give more details on an
efficient implementation of this algorithm, in particular on critical-pair
criteria specific to lattice ideal computations. Finally, we conclude the paper
with a computational comparison of our implementation of the Project-and-Lift
algorithm in 4ti2 with algorithms for lattice ideal computations implemented in
CoCoA and Singular. Our algorithm outperforms the other algorithms in every
single instance we have tried.
\end{abstract}
\end{frontmatter}

\section{Introduction}
In this article, we present a new algorithm for computing a
\emph{generating set} of a \emph{lattice ideal}
\[
I(\Lattice):= \langle x^{u^+} - x^{u^-}: u \in \Lattice \rangle
\subseteq k[x_1,...,x_n],
\]
where $k$ is a field, $\Lattice$ is a sub-lattice of $\Zset^n$, and
\[x^{u^+} - x^{u^-} := x_1^{u^+_1}x_2^{u^+_2} \cdot\cdot\cdot x_n^{u^+_n}
    - x_1^{u^-_1}x_2^{u^-_2} \cdot\cdot\cdot x_n^{u^-_n}\] where
$u^+_i = \max\{u_i,0\}$ and $u^-_i = \max\{-u_i,0\}$.
We assume that $\Lattice \cap \Nset^n = \{0\}$.
Generating set and Gr\"obner basis computations for general ideals are usually
very time consuming. Fortunately, in the special setting of lattice ideals,
many improvements are possible.  Two interesting areas of application of lattice
ideals are algebraic statistics and integer programming.

In general, a generating set of $I(\Lattice)$ is \emph{not} readily available.
For a basis $S$ of the lattice $\Lattice$ over $\Zset^n$,
the ideal $J(S) :=\langle x^{u^+}-x^{u^-}:u\in S\rangle$
satisfies $J(S) \subseteq I(\Lattice)$, but usually one may
not expect that $J(S)=I(\Lattice)$.  Also, when computing a Gr\"obner basis of a
lattice ideal, computational experiments show that when
$\Lattice \cap \Nset^n = \{0\}$, the computation of a generating set
usually takes much longer than computing the Gr\"obner basis from
the generating set.

Generating sets of lattice ideals and Gr\"obner bases of lattice ideals have
corresponding geometric concepts
(see \citet{Thomas:95, Urbaniak+Weismantel+Ziegler:97, Weismantel:98}),
which we call \emph{generating sets of lattices} and
\emph{Gr\"obner bases of lattices} respectively
(see Section \ref{Section: Generating sets and Groebner bases}).
These concepts are related as follows:
if a set $S \subseteq \Lattice$ is a generating set of $\Lattice$ or
a Gr\"obner basis of $\Lattice$ with respect to a term order $\succ$,
then $G := \{x^{u^+}-x^{u^-}: u \in S\}$ is respectively a generating set of
$I(\Lattice)$ or a Gr\"obner basis of $I(\Lattice)$ with
respect to $\succ$; and also,
if a set of monic binomials $G$ is a generating set of $I(\Lattice)$ or a
Gr\"obner basis of $I(\Lattice)$ with respect to a term order $\succ$,
then $S := \{\alpha - \beta: x^\alpha - x^\beta \in G\}$ is respectively
a generating set of $\Lattice$ or a Gr\"obner basis of $\Lattice$ with respect
to $\succ$.
Note that we use the same order $\succ$ for monomials and vectors via the
relation $x^\alpha \succ x^\beta$ if and only if $\alpha \succ \beta$
for $\alpha, \beta \in \Nset^n$. Also, note that any minimal reduced Gr\"obner
basis of a lattice ideal is a set of monic binomials.

In this paper, we have chosen to present existing
theory and the new algorithm only in a geometric framework following the
approach in \citet{Thomas:95}, \citet{Urbaniak+Weismantel+Ziegler:97},
and \citet{Weismantel:98}, since for lattice ideals, we prefer the geometric
approach to the algebraic one.
However, note that for every
geometric concept presented, there exists an equivalent
algebraic notion, although we do not present it here.
Also, we have tried to make this paper reasonably self contained, so for
completeness, we present geometric proofs of existing results where pertinent.

Recently, there has been renewed interest in toric ideal computations because of
applications in algebraic statistics.
Here, we are interested in Markov bases, which are used in a Monte-Carlo
Markov-Chain (MCMC) process to test validity of statistical models via sampling.
Diaconis and Sturmfels \citep{Diaconis+Sturmfels:98} showed 
that a (preferably minimal) generating set of a lattice
\[\Lattice_A := \{u\in\Zset^n:Au=0\}\] for some matrix
$A \in \Zset^{d \times n}$ is a Markov basis.
Note that $I(\Lattice_A)$ is a toric ideal.
However, at that time, no effective implementation of an algorithm to compute
generating sets of toric ideals was available that could deal with moderate size
problems in $50-100$ variables. This situation has changed by now:
several such implementations are available. Using \FourTiTwo\ \citep{4ti2},
Eriksson even reports, in \citet{Eriksson:04}, on successful computations
of Gr\"obner bases and Markov bases of toric ideals in $2,048$
variables. His problems arise from phylogenetic trees in
computational biology.

In integer programming, test sets of integer programs correspond to
Gr\"obner bases of lattices (or lattice ideals).
Test sets were introduced in \citet{Graver:75}.
Consider the general linear integer program
\[
\min\{c z_\sigmabar:Az=b, z_\sigmabar\geq 0, z\in\Zset^n\}
\]
where $c \in \Qset^{|\sigmabar|}$, $A \in \Zset^{d \times n}$, $b \in \Zset^d$,
$\sigma\subseteq\{1,...,n\}$, $\sigmabar:=\{1,...,n\}\setminus\sigma$,
and where $z_\sigmabar$ is the set of variables indexed by $\sigmabar$.
Any integer program that has an optimal solution can be written in this
form \citep{Conti+Traverso:91}.
By projecting onto the $\sigmabar$ variables, we can rewrite these
integer programs in the equivalent and more convenient form
\[
IP^\sigma_{A,c,b}= \min\{c z: A_\sigmabar z
\equiv b \pmod{A_\sigma\Zset}, z\in\Nset^{|\sigmabar|}\},
\]
where $A_\sigma$ and $A_\sigmabar$ are the sub-matrices of $A$ whose
columns are indexed by $\sigma$ and $\sigmabar$ respectively, and
$A_\sigma\Zset:=\{A_\sigma z: z\in \Zset^{|\sigma|}\}$.
In the special case where $\sigma=\emptyset$,
we set $A_\sigma\Zset:=\{0\}$, and the problem
$IP^\emptyset_{A,c,b}$ simplifies to
$IP_{A,c,b} := \min\{c z: Az=b, z\in\Nset^n\}$. Note that group relaxations and
extended group relaxations of $IP_{A,c,b}$ are also of the form
$IP^\sigma_{A,\bar{c},b}$ for some cost vector $\bar{c} \in \Qset^{|\sigmabar|}$
\citep{Hosten+Thomas}.
Without loss of generality, we assume that $c$ is \emph{generic} meaning that
$IP^\sigma_{A,c,b}$ has a unique optimal solution for every feasible
$b \in \Zset^d$.
We can always easily perturb a given $c$ so that it is generic.

A set $T\subseteq\Zset^{|\sigmabar|}$ is
called a \emph{test set} for $IP^\sigma_{A,c,b}$ if $T$ contains an
improving direction $t$ for every non-optimal feasible solution
$z \in \Nset^{|\sigmabar|}$ of $IP^\sigma_{A,c,b}$;
that is, $z-t$ is also feasible and $c (z-t)<c z$.
Clearly, $z-t$ being feasible implies that $t$ is an element of the lattice
\[
\Lattice^\sigma_A:= \{u \in \Zset^{|\sigmabar|} : A_\sigmabar u \equiv 0
\pmod{A_\sigma\Zset}\}.
\]
Moreover, a set $T\subseteq\Lattice^\sigma_A$ is called a test set for
$IP^\sigma_{A,c}$ the class of integer programs $IP^\sigma_{A,c,b}$
for all $b\in\Zset^d$, if $T$ is a test set for every integer program in
$IP^\sigma_{A,c}$.
Graver showed that there exist \emph{finite} sets $T$ that are test
sets for $IP_{A,c}$ ($\sigma = \emptyset$). In fact, his sets also constitute
finite test sets for $IP^\sigma_{A,c}$ for arbitrary $\sigma$. Having a finite
test set available, an optimal solution of $IP^\sigma_{A,c,b}$ can be found
by iteratively improving any given non-optimal solution of
$IP^\sigma_{A,c,b}$.

In \citet{Conti+Traverso:91} and \citet{Sturmfels+Weismantel+Ziegler:95},
it was shown that given a generic cost vector $c$ and a term order
$\succ$ where $c$ and $\succ$ are compatible,
a set $S \subseteq \Lattice^\sigma_A$ is a Gr\"obner basis of
$\Lattice^\sigma_A$ with respect to a $\succ$ if and only if $S$ is a test set
for $IP^\sigma_{A,c}$ where $c$ and $\succ$ are compatible
if $c \alpha > c \beta$ implies $\alpha \succ \beta$ for all
$\alpha,\beta \in \Nset^n$.
A compatible term ordering $\succ$ exists for every generic $c$, and a compatible
generic $c$ exists for every term ordering $\succ$.
Additionally, any lattice $\Lattice$ can be written in the form
$\Lattice^\sigma_A$ for some matrix $A \in \Zset^{n \times d}$ and some index
set $\sigma \subseteq \{1,...,n\}$, and so,
Gr\"obner bases of lattices and test sets of integer programs really are
equivalent concepts.

We define generating sets and Gr\"obner bases of lattices,
in Section \ref{Section: Generating sets and Groebner bases},
in a geometric context and present the completion procedure
\citep{Buchberger:87,Buchberger:85,Cox+Little+OShea:92}, which is the main
building block for the algorithms for computing generating sets.

In Section \ref{Section: Computing a generating set},
we present the two main existing algorithms for computing generating sets:
the algorithm of Hosten and Sturmfels in \citet{Hosten+Sturmfels}, which we call
the ``Saturation'' algorithm;
and the algorithm of Bigatti, LaScala, and Robbiano in
\citet{Bigatti+Lascala+Robbiano:99}, which we call the ``Lift-and-Project''
algorithm.
We also describe our new algorithm for computing generating sets:
the ``Project-and-Lift'' algorithm.
The Saturation Algorithm is based upon the result that
$ I(\Lattice)=(\ldots ((J(S):x_1^\infty):x_2^\infty)\ldots):x_n^\infty $
where $S$ is a lattice basis of $\Lattice$ and $J(S)$ is defined as above.
Using this result, we can compute a generating set of $I(\Lattice)$
from $S$ via a sequence of \emph{saturation} steps where each individual
saturation step is performed via the completion procedure.
The Lift-and-Project Algorithm is based upon the related result that
$I(\Lattice)=J(S):(x_1\cdot x_2 \cdot\ldots\cdot x_n)^{\infty}$. Here, a
generating set is computed via the completion procedure using an additional
variable. The Project-and-Lift Algorithm is strongly related to the
Saturation Algorithm; however, the computational speed-up is enormous as will be
seen in Section \ref{Section: Computational experience}.
In contrast with the Saturation Algorithm, which performs saturation
steps in the original space of the lattice $\Lattice$, the Project-and-Lift
Algorithm performs saturation steps in projected subspaces of $\Lattice$ and
then lifts the result back into the original space.

We are mainly interested in computing a generating set of $\Lattice$ where
$\Lattice \cap \Nset^n = \{0\}$.  However in Section \ref{Section: What if},
we address the question of how to compute a generating set $\Lattice$
if $\Lattice \cap \Nset^n \ne \{0\}$.
We demonstrate that the above methods for the case where
$\Lattice \cap \Nset^n = \{0\}$
can be extended to this more general case; it happens to be more
straight-forward in some ways.

The completion procedure as it is presented in Section
\ref{Section: Generating sets and Groebner bases} is not very efficient.
In Section \ref{Section: Speeding-up the Completion Procedure}, we show
how to increase the efficiency of the completion procedure.
All the results in this section, which we present in a geometric framework,
have corresponding results in an algebraic context.
This section is rather technical and no other section depends upon it,
so it may be skipped on first reading.

In Section \ref{Section: The 4x4x4-challenge}, we give the solution
of a computational challenge posed by Seth Sullivant to compute the
Markov basis of $4\times 4\times 4$ tables with $2$-marginals, a problem
involving $64$ variables. We solved this with the help of the new algorithm.
Our computations led to $148,968$ elements in the
minimal generating set of $I(A)$ which fall into $15$ equivalence
classes with respect to the underlying symmetry group $S_4\times
S_4\times S_4\times S_3$.

In Section \ref{Section: Computational experience}, we compare the
performance of the implementation of the Project-and-Lift
algorithm in \FourTiTwo\ v.1.2 \citep{4ti2} with the implementation
of the Saturation algorithm and the Lift-and-Project algorithm in Singular
v3.0.0 \citep{Singular}  and in CoCoA 4.2 \citep{CoCoA}.  The
Project-and-Lift algorithm is significantly faster than the other
algorithms.

\section{Generating sets and Gr\"obner bases}
\label{Section: Generating sets and Groebner bases}
Given a lattice $\Lattice \subseteq \Zset^n$, and a vector $b \in \Zset^n$,
we define
\[\Feasible_{\Lattice,b}:=\{x: x\equiv b\pmod{\Lattice}, x\in\Nset^n\}.\]
For $S\subseteq\Lattice$, we define $\Graph(\Feasible_{\Lattice,b},S)$ to be the
\emph{undirected} graph with nodes $\Feasible_{\Lattice,b}$ and edges $(x,y)$
if $x-y\in S$ or $y-x\in S$ for $x,y \in \Feasible_{\Lattice,b}$.
\begin{defn} \label{Definition: Generating set}
A set $S\subseteq\Lattice$ a \textbf{generating set of $\Lattice$}
if the graph $\Graph(\Feasible_{\Lattice,b},S)$ is connected for every
$b\in\Zset^n$.
\end{defn}
We remind the reader that connectedness of
$\Graph(\Feasible_{\Lattice,b},S)$ simply states that between
each pair $x,y\in \Feasible_{\Lattice,b}$, there exists a path
from $x$ to $y$ in $\Graph(\Feasible_{\Lattice,b},S)$.
Note the difference between a \emph{generating set} of a lattice and
a \emph{spanning set} of a lattice: 
a \textbf{spanning set of $\Lattice$} is any set $S \subseteq \Lattice$
such that any point in $\Lattice$ can be represented as an linear integer
combination of the vectors in $S$.
A generating set of $\Lattice$ is a spanning
set of $\Lattice$, but the converse is not necessarily true.

Recall that for any lattice $\Lattice$, we have 
$\Lattice = \Lattice^\sigma_A$ for some some matrix $A \in \Zset^{n \times d}$ and
some index set $\sigma \subseteq \{1,...,n\}$.
Hence, \[\Feasible_{\Lattice,b} = \Feasible^\sigma_{A,\bar{b}} :=
\{x\in\Nset^{|\sigmabar|}: A_\sigmabar x \equiv \bar{b} \pmod{A_\sigma\Zset}\}\]
for all $b\in\Zset^n$ and all $\bar{b} \in \Zset^d$ where $\bar{b} = A_\sigmabar b$.
So, $\Feasible_{\Lattice,b}$ and $\Feasible^\sigma_{A,\bar{b}}$ are 
dual representations of feasible sets.

\begin{exmp}
Let $S := $\emph{\{(1,-1,-1,-3,-1,2),(1,0,2,-2,-2,1)\}}, and
let $\Lattice \subseteq \Zset^6$ be the lattice spanned by $S$.
So, by definition, $S$ is a spanning set of $\Lattice$, but $S$ is not a
generating set of $\Lattice$.
Observe that $\Lattice = \Lattice_A$ where
\[
A = (\tilde{A},I),\;
\tilde{A} = 
\begin{pmatrix}
 -2 & -3 \\
 +2 & -1 \\
 +2 & +1 \\
 -1 & +1  
\end{pmatrix}, \text{ and }
I = 
\begin{pmatrix}
 1 & 0 & 0 & 0 \\
 0 & 1 & 0 & 0 \\
 0 & 0 & 1 & 0 \\
 0 & 0 & 0 & 1  
\end{pmatrix}.
\]
So, for every $b \in \Zset^6$,
$\Feasible_{\Lattice,b} = \Feasible_{A,\bar{b}} =
\{(x,s): \tilde{A}x + Is = \bar{b}, x\in\Nset^2,s\in\Nset^4\}$ where
$\bar{b} = Ab \in \Zset^4$.
Hence, the projection of $\Feasible_{\Lattice,b}$ onto the
$(x_1,x_2)$-plane is the set of integer points in the polyhedron
$\{x \in \Rset^n_+: Ax \le \bar{b}\}$,
and the $s$ variables are the slack variables.
Consider $b := (2,2,4,2,4,1)$;
then, $\Feasible_{\Lattice,b} = \Feasible_{A,\bar{b}}$ where
$\bar{b} = Ab = $\emph{(-6,4,10,1)}
(see Figure \ref{Figure: Connected Graph}a).

\input{Figure_Connected_Graph}

The graph
$\Graph(\Feasible_{\Lattice,b},S)$ is not connected
because the point $(3,4,12,2,0,0) \in \Feasible_{\Lattice,b}$ is disconnected
(see Figure \ref{Figure: Connected Graph}b). 
Let $S' := S \cup $\emph{\{(1,1,5,-1,-3,0)\}}. The graph of
$\Graph(\Feasible_{\Lattice,b},S')$ is now connected
(see Figure \ref{Figure: Connected Graph}c); however, $S'$ is still not a
generating set of $\Lattice$ since we have
$\Feasible_{\Lattice,b'}=\{(0,0,0,0,1,1),(0,1,3,1,0,0)\}$ 
for $b' := (0,0,0,0,1,1)$, and the graph
$\Graph(\Feasible_{\Lattice,b'},S')$ is disconnected; since there are only two
feasible points in $\Feasible_{\Lattice,b'}$, the vector between them
\emph{(0,1,3,1,-1,-1)} must be in any generating of $\Lattice$.
Finally, the set $S'' := S' \cup $\emph{\{(0,1,3,1,-1,-1)\}}
is a generating set of
$\Lattice$.

\end{exmp}

For the definition of a Gr\"obner basis, we need a term ordering
$\succ$ for $\Lattice$. We call $\succ$ a {\bf term ordering} for $\Lattice$
if $\succ$ is a total well-ordering on $\Feasible_{\Lattice,b}$ for every
$b \in \Zset^n$ and $\succ$ is an additive ordering meaning that
for all $b \in \Zset^n$ and for all $x,y \in \Feasible_{\Lattice,b}$, if 
$x \succ y$, then $x+\gamma\succ y+\gamma$ for every $\gamma\in\Nset^n$
(note that $x+\gamma, y+\gamma \in \Feasible_{\Lattice,b+\gamma}$).
We also need the notion of a decreasing path: a path $(x^0,\ldots,x^k)$ in
$\Graph(\Feasible_{\Lattice,b},G)$ is {\bf $\succ$-decreasing}
if $x^i\succ x^{i+1}$ for $i=0,\ldots,k-1$.
We define $\Lattice_\succ := \{u \in \Lattice: u^+ \succ u^-\}$.

\begin{defn}
A set $G\subseteq\Lattice_\succ$ is a {\bf $\succ$-Gr\"obner basis} of
$\Lattice$ if for every $x\in\Nset^n$ there
exists a $\succ$-decreasing path in $\Graph(\Feasible_{\Lattice,x},G)$
from $x$ to the unique $\succ$-minimal element in
$\Feasible_{\Lattice,x}$.
\end{defn}

If $G\subseteq\Lattice_\succ$ is a $\succ$-Gr\"obner basis, then $G$ is a
generating set of $\Lattice$ since given $x,y \in
\Graph(\Feasible_{\Lattice,b},G)$ for some $b \in \Zset^n$, there exists a
$\succ$-decreasing path from $x$ to the unique
$\succ$-minimal element in $\Feasible_{\Lattice,b}$ and from $y$ to the same
element, and thus, $x$ and $y$ are connected in
$\Graph(\Feasible_{\Lattice,b},G)$.
Also, $G\subseteq\Lattice_\succ$ is a Gr\"obner basis if and only if for every
$x\in\Nset^n$, $x$ is either the unique $\succ$-minimal element in
$\Feasible_{\Lattice,b}$ or there exists a vector $u \in G$ such that $x-u \in
\Feasible_{\Lattice,b}$ and $x \succ x-u$;
consequently, a Gr\"obner basis $G$ is a test set for $IP^\sigma_{A,c}$ 
where $\Lattice^\sigma_A = \Lattice$ if $c$ and $\succ$ are compatible.

The defining property of a Gr\"obner basis is very strong, so we
redefine it in terms of reduction paths. A path
$(x^0,\ldots,x^k)$ in $\Graph(\Feasible_{\Lattice,b},G)$
is a {\bf $\succ$-reduction path} if for \textbf{no}
$i\in\{1,\ldots,k-1\}$, we have $x^i\succ x^0$ and $x^i\succ x^k$.
For example, see Figure \ref{Figure: Reduction Path}.

\input{Figure_Reduction_Path}

\begin{lem} \label{Lemma: G is Groebner basis if there exist
reduction paths for all x and y} A set $G\subseteq\Lattice_\succ$
is a $\succ$-Gr\"obner basis of $\Lattice$ if
and only if for each $b\in\Zset^n$ and for each pair $x,y\in
\Feasible_{\Lattice,b}$, there exists a $\succ$-reduction path in
$\Graph(\Feasible_{\Lattice,b},G)$ between $x$ and $y$.
\end{lem}
\boproof
If $\Graph(\Feasible_{\Lattice,b},G)$ contains $\succ$-decreasing paths from
$x,y\in\Feasible_{\Lattice,b}$ to the unique $\succ$-minimal element
in $\Feasible_{\Lattice,b}$, then joining the two paths (and
removing cycles if necessary) forms a $\succ$-reduction path between $x$ and
$y$.

For the other direction, we assume that there is a $\succ$-reduction path
between each pair $x,y\in\Feasible_{\Lattice,b}$. Denote by $x^*$
the unique $\succ$-minimal element in $\Feasible_{\Lattice,b}$;
thus, every $x\in \Feasible_{\Lattice,b}$ is connected to $x^*$ by a
$\succ$-reduction path.
In particular, by the definition of a $\succ$-reduction
path, if $x\neq x^*$, then the first node $x^1\neq x$ in this path
must satisfy $x \succ x^1$. Repeating this argument iteratively with
$x^1$ instead of $x$, we get a $\succ$-decreasing path from $x$ to $x^*$.
This follows from the fact that $\succ$ is a term ordering, which
implies that every $\succ$-decreasing path must be finite. However, the only
node from which the $\succ$-decreasing path cannot be lengthened is $x^*$.
\eoproof

Checking for a given $G\subseteq\Lattice_\succ$ whether there
exists a $\succ$-reduction path in $\Graph(\Feasible_{\Lattice,b},G)$
for every $b\in\Zset^n$ and for each pair $x,y\in \Feasible_{\Lattice,b}$ involves
infinitely many situations that need to be checked. In fact, far fewer checks
are needed: we only need to check for a $\succ$-reduction path from $x$ to $y$
if there exists a \textbf{$\succ$-critical path} from $x$ to $y$.

\begin{defn}
Given $G \subseteq \Lattice_\succ$ and $b \in \Zset^n$, a path $(x,z,y)$
in $\Graph(\Feasible_{\Lattice,b},G)$
is a {\bf $\succ$-critical path} if $z\succ x$ and $z\succ y$.
\end{defn}
If $(x,z,y)$ is a $\succ$-critical path in $\Graph(\Feasible_{\Lattice,b},G)$,
then $x+u=z=y+v$ for some pair $u,v \in G$, in which case, we call $(x,z,y)$ a
$\succ$-critical path for $(u,v)$ (see Figure \ref{Figure: uv Path}).

\input{Figure_uv_Path}

The following lemma will be a crucial ingredient in the correctness
proofs of the algorithms presented in Section \ref{Section:
Computing a generating set}. It will guarantee correctness of the
algorithm under consideration, since the necessary reduction paths
have been constructed during the run of the algorithm. In the next lemma,
we cannot assume that $G$ is a generating set of $\Lattice$, since often
this is what we are trying to construct.

\begin{lem}
\label{Lemma: Test for reduction paths of critical paths is sufficient for
each fiber} Let $x,y\in \Feasible_{\Lattice,b}$ for some $b\in\Zset^n$, and let
$G\subseteq\Lattice_\succ$ where there is a path between
$x$ and $y$ in $\Graph(\Feasible_{\Lattice,b},G)$. If there
exists a $\succ$-reduction path between $x'$ and $y'$ for every $\succ$-critical
path $(x',z',y')$ in $\Graph(\Feasible_{\Lattice,b},G)$, then there exists a
$\succ$-reduction path between $x$ and $y$ in
$\Graph(\Feasible_{\Lattice,b},G)$.
\end{lem}
\boproof Assume on the contrary that no such $\succ$-reduction path
exists from $x$ to $y$. Among all paths $(x=x^0,\ldots,x^k=y)$ in
$\Graph(\Feasible_{\Lattice,b},G)$ choose one such that
$\max\limits_\succ\{x^0,\ldots,x^k\}$ is minimal. Such a minimal path
exists since $\succ$ is a term ordering. Let $j \in \{0,\dots,k\}$ where $x^j$
attains this maximum.

By assumption, $(x^0,\dots,x^k)$ is not a $\succ$-reduction path, and thus,
$x^j\succ x^0$ and $x^j\succ x^k$, and since $x^j$ is maximal,
we have $x^j\succ x^{j-1}$ and $x^j\succ x^{j+1}$. Let
$u=x^j-x^{j-1}$ and $v=x^j-x^{j+1}$. Then $(x^{j-1},x^j,x^{j+1})$
forms a $\succ$-critical path.
Consequently, we can replace the path $(x^{j-1},x^j,x^{j+1})$ with the
$\succ$-reduction path $(x^{j-1}=\bar{x}^0,\ldots,\bar{x}^s=x^{j+1})$ in the
path $(x^0,\ldots,x^k)$ and obtain a new path between $x$ and $y$ with the
property that the $\succ$-maximum of the intermediate nodes is strictly less
than $x^j=\max\limits_\succ\{x^1,\ldots,x^{k-1}\}$ (see Figure
\ref{Figure: Replacing Reduction Path}). This contradiction proves
our claim. \eoproof

\input{Figure_Replacing_Reduction_Path}

The following corollary is a straight-forward consequence of Lemma
\ref{Lemma: Test for reduction paths of critical paths is sufficient
for each fiber}, but nonetheless, it is worthwhile stating explicitly.
\begin{cor}
\label{Corollary: Test for reduction paths of critical paths is sufficient
for each fiber}
Let $G\subseteq\Lattice_\succ$.
If for all $b'\in\Zset^n$ and
for every $\succ$-critical path $(x',z',y')$ in
$\Graph(\Feasible_{\Lattice,b'},G)$,
there exists a $\succ$-reduction path between $x'$ and $y'$,
then for all $b \in \Zset^n$ and for all $x,y \in \Feasible_{\Lattice,b}$
where $x$ and $y$ are connected in $\Graph(\Feasible_{\Lattice,b},G)$,
there exists a $\succ$-reduction path between $x$ and $y$ in
${\Graph}(\Feasible_{\Lattice,b},G)$.
\end{cor}

Combining Corollary
\ref{Corollary: Test for reduction paths of critical paths is sufficient
for each fiber}
with Lemma \ref{Lemma: G is Groebner basis if there exist reduction
paths for all x and y}, we arrive at the following result for Gr\"obner bases.

\begin{cor}
\label{Corollary: Test for reduction paths of critical paths is
sufficient} A set $G\subseteq\Lattice_\succ$ is a $\succ$-Gr\"obner
basis of $\Lattice$ if and only if $G$ is a
generating set of $\Lattice$ and if for all $b\in\Zset^n$ and for every
$\succ$-critical path $(x,z,y)$ in $\Graph(\Feasible_{\Lattice,b},G)$, there
exists a $\succ$-reduction path between $x$ and $y$ in
$\Graph(\Feasible_{\Lattice,b},G)$.
\end{cor}

In Corollary
\ref{Corollary: Test for reduction paths of critical paths is sufficient
for each fiber} and Corollary
\ref{Corollary: Test for reduction paths of critical paths is
sufficient}, it is not necessary to check for a $\succ$-reduction path from $x$
to $y$ for every $\succ$-critical path $(x,y,z)$
in $\Graph(\Feasible_{\Lattice,b},G)$
for all $b\in\Zset^n$. Consider the case where there exists another
$\succ$-critical path $(x',y',z')$ in $\Graph(\Feasible_{\Lattice,b'},G)$
for some $b'\in\Zset^n$ such that $(x,y,z)=(x'+\gamma,y'+\gamma,z'+\gamma)$
for some $\gamma \in \Nset^n$. Then, a $\succ$-reduction path from
$x'$ to $y'$ in $\Graph(\Feasible_{\Lattice,b'},G)$ translates by $\gamma$
to a $\succ$-reduction path from $x$ to $y$ in
$\Graph(\Feasible_{\Lattice,b},G)$.
Thus, we only need to check for a $\succ$-reduction path from $x'$ to $y'$.

A $\succ$-critical path $(x,y,z)$ is {\bf minimal} if there does not exist
another $\succ$-critical path $(x',y',z')$ such that
$(x,y,z)=(x'+\gamma,y'+\gamma,z'+\gamma)$
for some $\gamma \in \Nset^n$ where $\gamma \ne 0$, or equivalently,
$\min\{x_i,y_i,z_i\}=0$ for all $i = 1,\dots,n$.
Consequently, if there exists a $\succ$-reduction path between $x$ and $y$ for
all minimal $\succ$-critical paths $(x,y,z)$, then there exists a
$\succ$-reduction path between $x'$ and $y'$ for all $\succ$-critical paths
$(x',y',z')$.
Also, for each pair of vectors $u,v \in \Lattice$, there
exists a unique minimal $\succ$-critical path $(x^{(u,v)},z^{(u,v)},y^{(u,v)})$
determined by $z^{(u,v)}:=\max\{u^+,v^+\}$ component-wise,
$x^{(u,v)}:=z^{(u,v)}-u$ and $y^{(u,v)}:=z^{(u,v)}-v$.
So, any other $\succ$-critical path for $(u,v)$ is of the form
$(x^{(u,v)}+\gamma,z^{(u,v)}+\gamma,y^{(u,v)}+\gamma)$ for some
$\gamma\in\Nset^n$. Using minimal $\succ$-critical paths, we can rewrite
Corollary \ref{Corollary: Test for reduction paths of critical paths is
sufficient for each fiber} and Corollary
\ref{Corollary: Test for reduction paths of critical paths is
sufficient}, so that we only need to check for a \emph{finite} number of
$\succ$-reduction paths.

\begin{lem}
\label{Lemma: Test for reduction paths of minimal critical paths is sufficient}
Let $G\subseteq\Lattice_\succ$.
If there exists a $\succ$-reduction path between $x^{(u,v)}$ and $y^{(u,v)}$ for
every pair $u,v \in G$, then
for all $b \in \Zset^n$ and for all $x,y \in \Feasible_{\Lattice,b}$
where $x$ and $y$ are connected in $\Graph(\Feasible_{\Lattice,b},G)$,
there exists a $\succ$-reduction path between $x$ and $y$ in
${\Graph}(\Feasible_{\Lattice,b},G)$
\end{lem}

\begin{cor}
\label{Corollary: Test for reduction paths of minimal critical paths is sufficient}
A set $G\subseteq\Lattice_\succ$ is a $\succ$-Gr\"obner
basis of $\Lattice$ if and only if $G$ is a
generating set of $\Lattice$ and for each pair $u,v\in G$, there
exists a $\succ$-reduction path between $x^{(u,v)}$ and
$y^{(u,v)}$ in $\Graph(\Feasible_{\Lattice,z^{(u,v)}},G)$.
\end{cor}

We now turn Lemma
\ref{Lemma: Test for reduction paths of minimal critical paths is sufficient}
into an algorithmic tool.
The following algorithm,
Algorithm \ref{Algorithm: Completion procedure} below,
called a completion procedure \citep{Buchberger:87}, guarantees
that if for a set $S\subseteq\Lattice$ the points $x$ and $y$ are connected
in $\Graph(\Feasible_{\Lattice,x},S)$, then there exists a
$\succ$-reduction path between $x$ and $y$ in
$\Graph(\Feasible_{\Lattice,x},G)$, where $G$ denotes the set
returned by the completion procedure.
Thus, if $S$ is a generating set of $\Lattice$, then Algorithm
\ref{Algorithm: Completion procedure} returns a set $G$ that is a
$\succ$-Gr\"obner basis of $\Lattice$ by Corollary
\ref{Corollary: Test for reduction paths of minimal critical paths is
sufficient}.

Given a set $S \subseteq \Lattice$, the completion procedure first sets $G:=S$
and then directs all vectors in $G$ according to $\succ$ such that
$G \subseteq \Lattice_\succ$.
Note that at this point $\Graph(\Feasible_{\Lattice,b},S) =
\Graph(\Feasible_{\Lattice,b},G)$ for all $b \in \Zset^n$.
The completion procedure then determines whether the set $G$ satisfies Lemma
\ref{Lemma: Test for reduction paths of minimal critical paths is sufficient};
in other words, it tries to find a reduction path from $x^{(u,v)}$ to
$y^{(u,v)}$ for every pair $u,v \in G$.
If $G$ satisfies Lemma
\ref{Lemma: Test for reduction paths of minimal critical paths is sufficient},
then we are done.
Otherwise, no $\succ$-reduction path was found for some $(u,v)$, in which case,
we add a vector to $G$ so that a $\succ$-reduction path exists, and then
again, test whether $G$ satisfies Lemma
\ref{Lemma: Test for reduction paths of minimal critical paths is sufficient},
and so on.

To check for a $\succ$-reduction path, using the ``Normal Form Algorithm'',
Algorithm \ref{algorithm: normal form} below,
we construct a maximal $\succ$-decreasing path
in $\Graph(\Feasible_{\Lattice,z^{(u,v)}},G)$ from
$x^{(u,v)}$ to some $x'$,
and a maximal $\succ$-decreasing path
in $\Graph(\Feasible_{\Lattice,z^{(u,v)}},G)$ from
$y^{(u,v)}$ to some $y'$.
If $x' = y'$,  then we have found a $\succ$-reduction path from
$x^{(u,v)}$ to $y^{(u,v)}$. Otherwise, we add the vector $r \in \Lattice_\succ$
to $G$ where $r := x'-y'$ if $x' \succ y'$, and $r := y'-x'$ otherwise, so
therefore, there is now a $\succ$-reduction path from
$x^{(u,v)}$ to $y^{(u,v)}$ in $\Graph(\Feasible_{\Lattice,z^{(u,v)}},G)$.
Note that before we add $r$ to $G$, since the paths from $x$ to $x'$ and from
$y$ to $y'$ are maximal, there does not
exist $u \in G$ such that $x' \ge u^+$ or $y' \ge u^+$. Therefore, there does
not exist $u \in G$ such that $r^+ \ge u^+$. This condition is needed to ensure
that the completion procedure terminates.

\begin{alg}{Normal Form Algorithm}
\label{algorithm: normal form}

\underline{Input:} a vector $x\in\Nset^n$ and a set $G\subseteq\Lattice_\succ$. 

\vspace{-0.35cm}

\underline{Output:} a vector $x'$ where there is a maximal $\succ$-decreasing path
from $x$ to $x'$ in $\Graph(\Feasible_{\Lattice,x},G)$.

\vspace{-0.35cm}

$x' := x$

\vspace{-0.35cm}

\underline{while} there is some $u\in G$ such that $u^+\leq x'$ \underline{do}

\vspace{-0.35cm}

\hspace{1.0cm} $x':=x'-u$

\vspace{-0.35cm}

\underline{return} $x'$
\end{alg}
We write $\NF(x,G)$ for the output of the Normal Form Algorithm.
\begin{alg}{Completion procedure}
\label{Algorithm: Completion procedure}

\underline{Input:} a term ordering $\succ$ and a set $S\subseteq\Lattice$.

\vspace{-0.35cm}

\underline{Output:} a set $G\subseteq\Lattice_\succ$ such that if $x,y$
are connected in $\Graph(\Feasible_{\Lattice,x},S)$, then there
exists a $\succ$-reduction path between $x$ and $y$ in
$\Graph(\Feasible_{\Lattice,x},G)$.

\vspace{-0.35cm}

\vspace{0.2cm} $G:=\{u:u^+\succ u^-, u\in S\}\cup\{-u:u^-\succ u^+, u\in S\}$

\vspace{-0.35cm}

$C:=\{(u,v): u,v \in G\}$

\vspace{-0.35cm}

\underline{while} $C\neq \emptyset $ \underline{do}

\vspace{-0.35cm}

\hspace{1.0cm}Select $(u,v) \in C$

\vspace{-0.35cm}

\hspace{1.0cm}$C:=C\setminus\{(u,v)\}$

\vspace{-0.35cm}

\hspace{1.0cm}$r:=\NF(x^{(u,v)},G)-\NF(y^{(u,v)},G)$

\vspace{-0.35cm}

\hspace{1.0cm}\underline{if} $r\neq 0$ \underline{then}

\vspace{-0.35cm}

\hspace{2.0cm}\underline{if} $r^- \succ r^+$ \underline{then} $r:=-r$

\vspace{-0.35cm}

\hspace{2.0cm}$C:=C\cup\{(r,s):s\in G\}$

\vspace{-0.35cm}

\hspace{2.0cm}$G:=G\cup \{r\}$

\vspace{-0.35cm}

\underline{return} $G$.
\end{alg}
We write $\CP(\succ,S)$ for the output of the Completion Procedure.

\begin{lem} \label{Lemma: One-by-one algorithm terminates and is
correct} Algorithm \ref{Algorithm: Completion procedure} terminates
and satisfies its specifications.
\end{lem}

\boproof
Let $(r^1,r^2,\dots)$ be the sequence of vectors $r$ that are added to
the set $G$ during the Algorithm \ref{Algorithm: Completion procedure}.
Since before we add $r$ to $G$, there does not exist $u \in G$
such that $r^+ \ge u^+$, the sequence satisfies ${r^i}^+\not\leq {r^j}^+$
whenever $i<j$. By the Gordan-Dickson Lemma (see for example
\citet{Cox+Little+OShea:92}), such a sequence must be finite and
thus, Algorithm \ref{Algorithm: Completion procedure} must
terminate.

When the algorithm terminates, the set $G$ must satisfy the property that
for each $u,v\in G$, there exists a $\succ$-reduction path from $x^{(u,v)}$ to
$y^{(u,v)}$, and therefore, by Lemma
\ref{Lemma: Test for reduction paths of minimal critical paths is sufficient},
there exists a $\succ$-reduction path between $x$ and $y$ in
${\Graph}(\Feasible_{\Lattice,b},G)$
for all $x,y \in \Feasible_{\Lattice,b}$ for all $b \in \Zset^n$
where $x$ and $y$ are connected in $\Graph(\Feasible_{\Lattice,b},G)$.
Moreover, by construction, $S \subseteq G \cup -G$, and therefore,
if $x$ and $y$ are connected in
$\Graph(\Feasible_{\Lattice,b},S)$, then $x$ and $y$ are connected in
$\Graph(\Feasible_{\Lattice,b},G)$.
\eoproof

Note that the completion procedure preserves connectivity:
given $x,y \in \Feasible_{\Lattice,b}$ for some $b \in \Zset^n$, if $x$ and $y$ are
connected in $\Graph(\Feasible_{\Lattice, b},S)$, then $x$ and $y$ are also
connected in $\Graph(\Feasible_{\Lattice, b},G)$.

\section{Computing a generating set}
\label{Section: Computing a generating set}

In this section, we finally present three algorithms to compute a
generating set of $\Lattice$: the ``Saturation'' algorithm
(Hosten and Sturmfels \citet{Hosten+Sturmfels}), 
the ``Lift-and-Project'' algorithm (Bigatti, LaScala, and Robbiano
\citet{Bigatti+Lascala+Robbiano:99}), and our \emph{new}
``Project-and-Lift'' algorithm.
Each algorithm produces a generating set of $\Lattice$ that is not
necessarily minimal, and so,
once a generating set of $\Lattice$ is known, a minimal
generating set of $\Lattice$ can be computed by a single Gr\"obner
basis computation (see \citet{Caboara+Kreuzer+Robbiano:03} for more
details).
The fundamental idea behind all three algorithms is essentially the same, and
the main algorithmic building block of the algorithms is the completion
procedure.

\subsection{The ``Saturation'' algorithm}
\label{Subsection: The Saturation-algorithm}
Let $x,y \in \Feasible_{\Lattice, b}$ for some $b\in\Zset^n$,
and let $S \subseteq \Lattice$.
Observe that if $x$ and $y$ are connected in
$\Graph(\Feasible_{\Lattice,b},S)$, then $x+\gamma$ and $y+\gamma$ are connected
in $\Graph(\Feasible_{\Lattice,b},S)$ for any $\gamma \in \Nset^n$ since we can
just translate any path from $x$ to $y$ in $\Graph(\Feasible_{\Lattice,b},S)$
by $\gamma$ giving a path from $x+\gamma$ to $y+\gamma$ in
$\Graph(\Feasible_{\Lattice,b+\gamma},S)$.
However,
it is not necessarily true that $x+\gamma$ and $y+\gamma$ are also connected
for any $\gamma \in \Zset^n$ ($\gamma$ may be negative)
where $x+\gamma \ge 0$ and $y+\gamma \ge 0$.

Given a set $S\subseteq\Lattice$, the Saturation algorithm constructs a set $T$
such that if $x$ and $y$ are connected in $\Graph(\Feasible_{\Lattice,b},S)$ for
some $b\in\Zset^n$, then $x+\gamma$ and $y+\gamma$ are connected
in $\Graph(\Feasible_{\Lattice,b},T)$
for any $\gamma \in \Zset^n$ where $x+\gamma \ge 0$ and $y+\gamma \ge 0$.
Importantly then, if $S$ spans $\Lattice$, then $T$ must be a generating set of
$\Lattice$. This follows since if $S$ spans $\Lattice$, then for all $b\in\Zset^n$
and for all $x,y \in \Feasible_{\Lattice, b}$, there must exist
a $\gamma \in \Nset^n$ such that $x+\gamma$ and $y+\gamma$ are connected in
$\Graph(\Feasible_{\Lattice,b},T)$, and hence, $x$ and $y$ must also be
connected in $\Graph(\Feasible_{\Lattice,b},T)$ from our assumption about $T$.

For convenience, we need some new notation.
Given $x,y \in \Nset^n$, we define $x \wedge y$ as the component-wise minimum of
$x$ and $y$ -- that is, $(x \wedge y)_i = \min\{x_i,y_i\}$ for all
$i=1,\dots,n$. Also, given $\sigma \subseteq \{1,\dots,n\}$, we define
$x \wedge_\sigma y$ as the component-wise minimum of $x$ and $y$ for the
$\sigma$ components and 0 otherwise -- that is,
$(x \wedge y)_i = \min\{x_i,y_i\}$ if $i\in \sigma$ and $(x \wedge y)_i = 0$
otherwise.
\begin{defn}
Let $\sigma \subseteq \{1,\dots,n\}$, and let $S,T \subseteq \Lattice$.
The set $T$ is \textbf{$\sigma$-saturated on $S$} if and only if
for all $b\in\Zset^n$ and for all $x, y \in \Feasible_{\Lattice, b}$ where
$x$ and $y$ are connected in $\Graph(\Feasible_{\Lattice,b},S)$, the points
$x-\gamma$ and $y-\gamma$ are also connected in
$\Graph(\Feasible_{\Lattice,b-\gamma},T)$ where
$\gamma = x \wedge_\sigma y$.
\end{defn}
So, when $T$ is $\sigma$-saturated on $S$, if $x$ and $y$
are connected in $\Graph(\Feasible_{\Lattice,b},S)$, then
$x+\gamma$ and $y+\gamma$ are connected in $\Graph(\Feasible_{\Lattice,b},T)$
for any $\gamma\in \Zset^n$ ($\gamma$ can be negative) where $x+\gamma\ge0$,
$y+\gamma\ge0$, and $\supp(\gamma) \subseteq \sigma$.
Saturation is thus concerned with the connectivity
of a set $T$ in relation to the connectivity of another set $S$.
Note that, by definition, a set $S \subseteq \Lattice$ is $\emptyset$-saturated
on itself.
Also observe that if $S$ spans $\Lattice$, then
$T \subseteq \Lattice$ is $\{1,\dots,n\}$-saturated on $S$ if and only if
$T$ is a generating set of $\Lattice$.

The fundamental idea behind the Saturation algorithm is given
$S,T \subseteq \Lattice$ where $T$ is $\sigma$-saturated on $S$ for some
$\sigma\subseteq \{1,\dots,n\}$,
we can compute a set $T'$ that is a $(\sigma\cup\{i\})$-saturated on $S$ for any
$i \in \sigmabar$.  Therefore, given a set $S \subseteq \Lattice$ that spans
$\Lattice$, starting from a set $T = S$, which is $\emptyset$-saturated on
$S$, if we do this repeatedly for each $i \in \{1,\dots,n\}$, we arrive at a
set $T' \subseteq \Lattice$ that is $\{1,\dots,n\}$-saturated on $S$ and,
therefore, a generating set of $\Lattice$.

The following two lemmas are fundamental to the saturation algorithm.
First, we extend the definition of reduction paths.
Given $\varphi \in \Qset^n$, a path $(x^0,\ldots,x^k)$ in
$\Graph(\Feasible_{\Lattice,b},G)$ is an \textbf{$\varphi$-reduction path}
if for \textbf{no} $j\in\{1,\ldots,k-1\}$,
we have $\varphi x^j > \varphi x^0$ and $\varphi x^j > \varphi x^k$.
Also, we define $e^i$ to be the $i$th unit vector and $\bar{e}^i = -e^i$.
So, given $b \in \Zset^n$, the path
$(x^0,\dots,x^k) \subseteq \Feasible_{\Lattice,b}$
is a $\bar{e}^i$-reduction path
if $x^j_i \ge x^0_i$ or $x^j_i \ge x^k_i$ for $j=1,\dots,k-1$.
\begin{lem}\label{Lemma: saturation}
Let $S,T \subseteq \Lattice$ and $i \in \{1,...,n\}$.
The set $T$ is $\{i\}$-saturated on $S$ if and
only if for all $b\in\Zset^n$ and for all $x, y \in \Feasible_{\Lattice, b}$ where 
$x$ and $y$ are connected in $\Graph(\Feasible_{\Lattice,b},S)$,
there exists a $\bar{e}^i$-reduction path from $x$ to $y$ in
$\Graph(\Feasible_{\Lattice,b},T)$.
\end{lem}
\boproof
Let $x,y \in \Feasible_{\Lattice, b}$ for some $b \in \Zset^n$ where $x$ and $y$
are connected in $\Graph(\Feasible_{\Lattice,b},S)$
and let $\gamma = x \wedge_{\{i\}} y$.

Assume $T$ is $\{i\}$-saturated on $S$, and so,
$x-\gamma$ and $y-\gamma$
are connected in $\Graph(\Feasible_{\Lattice,b-\gamma},T)$.
Let $(x-\gamma=x^0,\dots,x^k=y-\gamma)$ be a path from
$x-\gamma$ and $y-\gamma$ in $\Graph(\Feasible_{\Lattice,b-\gamma},T)$.
The path $(x=x^0+\gamma,\dots,x^k+\gamma=y)$ is a $\bar{e}^i$-reduction
path from $x$ to $y$ in $\Graph(\Feasible_{\Lattice,b},T)$.

Conversely, by assumption, there exists a $\bar{e}^i$-reduction path
$(x=x^0,\dots,x^k=y)$ in $\Graph(\Feasible_{\Lattice,b},T)$.
The path $(x-\gamma=x^0-\gamma,\dots,x^k-\gamma=y-\gamma)$ is thus a feasible 
path from $x-\gamma$ to $y-\gamma$ in $\Graph(\Feasible_{\Lattice,b-\gamma},T)$.
\eoproof

Given any vector $\varphi \in \Qset^n$ and a term order $\prec$, we define
the order $\prec_\varphi$ where $x \prec_\varphi y$ if $\varphi x < \varphi y$
or $\varphi x = \varphi y$ and $x \prec y$.
Since we assume $\Lattice \cap \Nset^n = \{0\}$, the order $\prec_{\bar{e}^i}$ is
thus a term ordering for $\Lattice$.
Importantly then, a $\prec_{\bar{e}^i}$-reduction path is also an
$\bar{e}^i$-reduction path.
Let $T=\CP(\prec_{\bar{e}^i},S)$.
Then, by the properties of the completion procedure, for all $b\in\Zset^n$ and
$x,y \in \Feasible_{\Lattice,b}$ where $x$ and $y$ are connected in
$\Graph(\Feasible_{\Lattice,b},S)$, there exists a $\prec_{\bar{e}^i}$-reduction
path from $x$ to $y$ in $\Graph(\Feasible_{\Lattice, b},T)$;
$T$ is therefore $\{i\}$-saturated on $S$.

Let $S,T \subseteq \Lattice$ and $T$ is $\sigma$-saturated on $S$ for some
$\sigma \subseteq \{1,\dots,n\}$.
Let $T'=\CP(\prec_{\bar{e}^i},T)$.
So, $T'$ is therefore $\{i\}$-saturated on $T$. For the saturation algorithm to
work, we need that $T'$ is also $(\sigma\cup\{i\})$-saturated on $S$, which
follows from Lemma \ref{Lemma: saturation union} below.

\begin{lem}\label{Lemma: saturation union}
Let $\sigma,\tau \subseteq \{1,\dots,n\}$ and
$S,T,U \subseteq \Lattice$.
If $U$ is $\sigma$-saturated on $S$, and $T$ is $\tau$-saturated on $U$,
then $T$ is $(\sigma\cup\tau)$-saturated on $S$.
\end{lem}
\boproof
Let $b \in \Zset^n$, and $x,y \in \Feasible_{\Lattice,b}$ where $x$ and $y$ are
connected in $\Graph(\Feasible_{\Lattice,b},S)$.
Let $\alpha = x \wedge_\sigma y$.
Since $T$ is $\sigma$-saturated on $S$, $x-\alpha$ and $y-\alpha$ are connected
in $\Graph(\Feasible_{\Lattice,b-\alpha},T)$.
Let $\beta  = x-\alpha \wedge_\tau y-\alpha$.
Then, since $U$ is $\tau$-saturated
on $T$, $x-\alpha-\beta$ and $y-\alpha-\beta$ are connected in
$\Graph(\Feasible_{\Lattice,b-\alpha-\beta},U)$.
Let $\gamma = \alpha+\beta$; then, $\gamma = x \wedge_{(\sigma \cup \tau)} y$.
Therefore, there is a path from $x-\gamma$ to $y-\gamma$
in $\Graph(\Feasible_{\Lattice,b-\gamma},U)$ as required.
\eoproof

We now arrive at the Saturation algorithm below.

\begin{alg}{Saturation algorithm}
\label{Saturation algorithm}

\underline{Input:} a spanning set $S$ of $\Lattice$.

\vspace{-0.35cm}

\underline{Output:} a generating set $G$ of $\Lattice$.

\vspace{-0.35cm}

$G:=S$

\vspace{-0.35cm}

$\sigma:=\emptyset$

\vspace{-0.35cm}

\underline{while} $\sigma\neq\{1,...,n\}$ \underline{do}

\vspace{-0.35cm}

\hspace{1.0cm}Select $i \in \sigmabar$

\vspace{-0.35cm}

\hspace{1.0cm}$G:=\CP(\prec_{\bar{e}^i},G)$

\vspace{-0.35cm}

\hspace{1.0cm}$\sigma:=\sigma\cup\{i\}$

\vspace{-0.35cm}

\underline{return} $G$.
\end{alg}

\begin{lem}
\label{Lemma: Saturation algorithm terminates and is correct}
Algorithm \ref{Saturation algorithm} terminates and satisfies its
specifications.
\end{lem}
\boproof
Algorithm \ref{Saturation algorithm} terminates, since
Algorithm \ref{Algorithm: Completion procedure} always terminates.
We show at the beginning of each iteration that $G$ is $\sigma$-saturated on
$S$, and so, at the end of the algorithm $G$ is $\{1,\dots,n\}$-saturated on
$S$;
therefore, $G$ is a generating set of $\Lattice$.
At the beginning of the first iteration, $G$ is $\sigma$-saturated on $S$ since
$\sigma = \emptyset$ and $G=S$. So, we can assume it is true for the current
iteration, and now, we show it is true for the next iteration.
Let $G':=\CP(\prec_{\bar{e}^i},G)$.
Then, by Lemma \ref{Lemma: saturation},
$G'$ is $\{i\}$-saturated on $G$,
and so, by Lemma \ref{Lemma: saturation union}, $G'$ is
$(\sigma\cup\{i\})$-saturated on $S$.  So, $G$ is $\sigma$-saturated on $S$ at
the beginning of the next iteration.
\eoproof

During the Saturation algorithm, we saturate $n$ times,
once for each $i \in \{1,\dots,n\}$.
However, as proven in \citet{Hosten+Shapiro}, it is in fact only necessary to
perform at most $\lfloor \frac{n}{2} \rfloor$ saturations.
Given $S,T \subseteq \Lattice$, we can show that there always
exists a $\sigma \subseteq \{1,\dots,n\}$ where
$|\sigma| \le \lfloor \frac{n}{2} \rfloor$ such that if $T$ is
$\sigma$-saturated on $S$, then $T$ is $\{1,\dots,n\}$-saturated on $S$.
The following two lemmas prove the result.

\begin{lem}\label{Lemma: support saturation}
Let $\sigma \subseteq \{1,\dots,n\}$, $S,T \subseteq \Lattice$
where $T$ is $\sigma$-saturated on $S$, and $u \in S$.
If $\supp(u^-) \subseteq \sigma$ or $\supp(u^+)\subseteq \sigma$,
then $T$ is $(\supp(u)\cup\sigma)$-saturated on $S$.
\end{lem}
\boproof
Assume that $\supp(u^-) \subseteq \sigma$.
Let $x,y \in \Feasible_{\Lattice,b}$ for some $b \in \Zset^n$
where $x$ and $y$ are connected in $\Graph(\Feasible_{\Lattice, b},S)$.
Let $\alpha = x \wedge_{\supp(u^+)} y$ and
$\beta = x-\alpha \wedge_\sigma y-\alpha$.
We must show that $x-\alpha-\beta$ and $y-\alpha-\beta$ are connected in 
$\Graph(\Feasible_{\Lattice, b-\alpha-\beta},T)$ since
$\alpha+\beta = x \wedge_{(\supp(u^+)\cup\sigma)} y$.
By translating the path from $x$ to $y$ by $\alpha$, we get a path from
$x-\alpha$ to $y-\alpha$ that is non-negative on all components except
$\supp(u^+)$.  This path can transformed into a path that is non-negative on all
components except $\supp(u^-)$ by adding $u$ to the start of the path as many
times as necessary and subtracting $u$ from the end of the path the same number
of times.
Therefore, $x-\alpha+\gamma$ and $y-\alpha+\gamma$ are connected in
$\Graph(\Feasible_{\Lattice, b-\alpha+\gamma},S)$
for some $\gamma \in \Nset^n$ where $\supp(\gamma) \subseteq \supp(u^-)
\subseteq \sigma$.  Observe that $\supp(\beta+\gamma) \subseteq \sigma$.
Thus, since $T$ is $\sigma$-saturated, $x-\alpha-\beta$ and $y-\alpha-\beta$ are
connected in $\Graph(\Feasible_{\Lattice, b-\alpha-\beta},T)$ as required.

The case where $T$ is $\supp(u^+) \subseteq \sigma$ is essentially the same as
above.
\eoproof

\begin{lem}\label{Lemma: upper bound on saturations} 
Let $S,T \subseteq \Lattice$.
There exists a $\sigma \subseteq \{1,\dots,n\}$
where $|\sigma| \le \lfloor \frac{n}{2} \rfloor$ such that if $T$ is
$\sigma$-saturated on $S$, then $T$ is $\{1,\dots,n\}$-saturated on $S$.
\end{lem}
\boproof
We show this by construction.
Without loss of generality, we assume that
$\Lattice$ is not contained in any of the linear subspaces
$\{x_i:x_i=0,x\in\Rset^n\}$ for $i=1,\ldots,n$; otherwise, we may simply
delete this component.

Let $\sigma = \emptyset$, $\tau = \emptyset$, and $U = \emptyset$. Repeat
the following steps until $\tau = \{1,\dots,n\}$.
\begin{enumerate}
\item Select $u \in S$ such that $\supp(u) \setminus \tau \ne \emptyset$.
\item If $|\supp(u^+) \setminus \tau| \ge |\supp(u^-) \setminus \tau|$,
then $\sigma := \sigma \cup \supp(u^-)$,
else $\sigma := \sigma \cup \supp(u^+)$.
\item Set $\tau := \tau \cup \supp(u)$, and set $U := U \cup \{u\}$.
\end{enumerate}
The procedure must terminate since during
each iteration we increase the size of $\tau$.
Note that, at termination, $U \subseteq S$,
$\bigcup_{u\in U}\supp(u)=\tau=\{1,\dots,n\}$,
and for all $u \in U$ either $\supp(u^+) \subseteq \sigma$ or
$\supp(u^-) \subseteq \sigma$.
Therefore, by applying Lemma \ref{Lemma: support saturation} recursively for
each $u \in U$, we have that if $T$ is $\sigma$-saturated
on $S$, then $T$ is $\{1,\dots,n\}$-saturated on $S$.
Lastly, since in each iteration we add at least twice as many components to
$\tau$ as to $\sigma$,
we conclude that at termination $|\sigma| \le \lfloor \frac{n}{2} \rfloor$.
\eoproof

\begin{exmp}
Consider again the set
$S :=$ \emph{\{(1,-1,-1,-3,-1,2),(1,0,2,-2,-2,1)\}}.
Let $\Lattice$ be the lattice spanned by $S$, and let $\sigma = \{1,6\}$.
Then, by Lemma \ref{Lemma: support saturation}, 
since $\supp(($\emph{1,-1,-1,-3,-1,2}$)^+) = \{1,6\}$ and
$\supp($\emph{(1,-1,-1,-3,-1,2)}$) = \{1,2,3,4,5,6\}$,
if a set $T \subseteq \Lattice$ is $\{1,6\}$-saturated on $S$, then $T$ is
$\{1,2,3,4,5,6\}$-saturated on $S$. So, to compute a generating set of
$\Lattice$, we only need to saturate on $\{1,6\}$.
The following table gives the values of $\sigma$, $i$, and $G$ at each stage of
the Saturation algorithm when constructing a set that is $\{1,6\}$-saturated on
$S$, and hence, a generating set of $\Lattice$.
\emph{\small{
\begin{center}
\begin{tabular}{|l|l|l|}
\hline
$\sigma$ & $i$ & $G := \CP(\prec_{\bar{e}^i},G)$ \\
\hline
\hline
$\emptyset$ & $1$ &  \{(-1,0,-2,2,2,-1) (-1,1,1,3,1,-2)
                          (-1,2,4,4,0,-3) (0,-1,-3,-1,1,1)\} \\
\hline
$\{1\}$     & $6$ &  \{(0,1,3,1,-1,-1), (-1,1,1,3,1,-2), (-1,0,-2,2,2,-1),
                       (-1,-1,-5,1,3,0), (1,2,8,0,-4,-1)\} \\
\hline
\end{tabular}
\end{center}
}}
Observe that after the first iteration, that $G$ is not a generating set of
$\Lattice$. The set $G$ does not contain the vector \emph{(-1,-1,-5,1,3,0)}, and
so, the graph $\Graph(\Feasible_{\Lattice,b},G)$ where $b = (0,0,0,1,3,0)$ is
disconnected.
Note that the final set $G$ is not a \emph{minimal} generating set of
$\Lattice$; the vector \emph{(1,2,8,0,-4,-1)} is not needed.
See \citet{Caboara+Kreuzer+Robbiano:03} for an algorithm to compute a minimal
generating set.
\end{exmp}

We now introduce the concept of a $\sigma$-generating set of $\Lattice$ for some
$\sigma \subseteq \{1,\dots,n\}$ --
a generalization of a generating set of $\Lattice$.
These new generating sets provide useful insights into saturation and 
the inspiration for the Project-and-Lift
algorithm as well as a point of reference to compare the two algorithms.

Firstly, we define
$\Feasible^\sigma_{\Lattice,b}:=\{z:z\equiv b\pmod{\Lattice},   
        z_\sigmabar\ge0,z\in\Zset^n\}$
where $\sigma \subseteq \{1,\dots,n\}$,
so we now allow the $\sigma$ components to be negative.
Given $S\subseteq\Lattice$, analogous to
$\Graph(\Feasible_{\Lattice,b},S)$, we define
$\Graph(\Feasible^\sigma_{\Lattice,b},S)$ to be the undirected graph with nodes
$\Feasible^\sigma_{\Lattice,b}$ and edges $(x,y)$ if $x-y\in S$ or $y-x \in S$.
Observe that a path in $\Graph(\Feasible^\sigma_{\Lattice,b},S)$ is
non-negative on the $\sigmabar$ components and may be negative
on the $\sigma$ components.
Analogous to a generating set of $\Lattice$, a set $S\subseteq\Lattice$ is a
$\sigma$-\textbf{generating set of $\Lattice$} if the graph
$\Graph(\Feasible^\sigma_{\Lattice,b},S)$ is connected for every $b\in\Zset^n$.
Note that $\emptyset$-generating sets are equivalent to generating sets
and $\{1,\dots,n\}$-generating sets are equivalent to spanning sets.
\begin{lem}\label{Lemma: sat to gen}
Let $\sigma \subseteq \{1,\dots,n\}$ and $S,T \subseteq \Lattice$ where $S$
spans $\Lattice$. If $T$ is $\sigma$-saturated on $S$, then $T$
is a $\sigmabar$-generating set of $\Lattice$.
\end{lem}
\boproof
Let $x,y \in \Feasible^\sigmabar_{\Lattice,b}$ for some $b \in \Zset^n$. We must show that
$x$ and $y$ are connected in $\Graph(\Feasible^\sigmabar_{\Lattice,b},S)$.
Since $S$ spans $\Lattice$, there must exist a
$\gamma \in \Nset^n$ such that $x+\gamma$ and $y+\gamma$ are connected in
$\Graph(\Feasible_{\Lattice,b+\gamma},S)$.
Let $\alpha,\beta \in \Nset^n$ where $\alpha+\beta=\gamma$,
$\supp(\alpha) \subseteq \sigma$, and $\supp(\beta) \subseteq \sigmabar$.
Since $T$ is $\sigma$-saturated on $S$ and $\supp(\alpha) \subseteq \sigma$,
the points $x+\beta=x+\gamma-\alpha$ and $y+\beta=y+\gamma-\alpha$ are connected
in $\Graph(\Feasible_{\Lattice,b+\beta},S)$.
Therefore, since $\supp(\beta) \subseteq \sigmabar$,
$x$ and $y$ are connected in $\Graph(\Feasible^\sigmabar_{\Lattice,b},S)$.
\eoproof

Interestingly, the converse of Lemma \ref{Lemma: sat to gen} is not true
in general: a $\sigmabar$-generating set is not necessarily a
$\sigma$-saturated set.
Let $S,T \subseteq \Lattice$ where $S$ spans $\Lattice$,
$\sigma \subseteq \{1,\dots,n\}$,
and let $x,y \in \Feasible_{\Lattice,b}$ for some $b\in\Zset^n$
where $x$ and $y$ are connected in $\Graph(\Feasible_{\Lattice,b},S)$.
If $T$ is a $\sigmabar$-generating set of $\Lattice$, then
$x-\gamma$ and $y-\gamma$ are connected in
$\Graph(\Feasible_{\Lattice^\sigmabar,b-\gamma},T)$ 
where $\gamma = x\wedge_\sigma y$.
In other words, there is a path from $x-\gamma$ and $y-\gamma$ that remains
non-negative on the $\sigma$ components but may be negative on the $\sigmabar$
components.
On the other hand,
if $T$ is $\sigma$-saturated on $S$,
then $x-\gamma$ and $y-\gamma$ are connected in
$\Graph(\Feasible_{\Lattice,b-\gamma},T)$
where again $\gamma = x\wedge_\sigma y$.
In other words, there is a path from $x-\gamma$ and $y-\gamma$ that remains
non-negative on \emph{all} the components.
So, while $\sigmabar$-generating sets, like $\sigma$-saturated sets,
ensure path non-negativity on the $\sigma$ components, they do not preserve
existing path non-negativity on the other $\sigmabar$ components like
$\sigma$-saturated sets do.
Indeed, $\sigmabar$-generating sets say nothing at all about the path
non-negativity of the $\sigmabar$ components. So, $\sigma$-saturation is a
stronger concept than $\sigmabar$-generation.

In the Project-and-Lift algorithm, we compute $\sigmabar$-generating sets
instead of $\sigma$-saturated sets. By doing so, we can effectively
\emph{ignore} the $\sigmabar$ components, and therefore, we compute smaller
intermediate sets, although we start and finish at the same point.

\subsection{The ``Project-and-Lift'' algorithm}

Given $\sigma \subseteq \{1,\dots,n\}$, we define the projective map
$\pi_\sigma: \Zset^n \mapsto \Zset^{|\sigmabar|}$ that projects a vector in $\Zset^n$
onto the $\sigmabar = \{1,\dots,n\}\setminus\sigma$ components.
For convenience, we write $\Lattice^\sigma$ where $\sigma\subseteq\{1,\ldots,n\}$
as the projection of $\Lattice$ onto the $\sigmabar$
components -- that is, $\Lattice^\sigma = \pi_\sigma(\Lattice)$.
Note that $\Lattice^\sigma$ is also a lattice.

The fundamental idea behind the Project-and-Lift algorithm is that using a set
$S \subseteq \Lattice^{\{i\}}$ that is a generating set of $\Lattice^{\{i\}}$
for some $i \in \{1,\dots,n\}$, we can compute a set
$S'\subseteq \Lattice^{\{i\}}$ such that $S'$ lifts
to a generating set of $\Lattice$.
So, for some $\sigma \subseteq \{1,\dots,n\}$, since $\Lattice^\sigma$ is also a
lattice, starting with a generating set of $\Lattice^\sigma$, we can compute a
generating of $\Lattice^{\sigma\backslash\{i\}}$ for some $i \in \sigma$.
So, by doing this repeatedly for every $i \in \sigma$,
we attain a generating set of $\Lattice$.

First, we extend the definition of Gr\"obner bases.
Given $\varphi \in \Qset^n$, recall that a path $(x^0,\ldots,x^k)$ in
$\Graph(\Feasible_{\Lattice,b},G)$ is a $\varphi$-reduction path
if for \textbf{no} $j\in\{1,\ldots,k-1\}$,
we have $\varphi x^j > \varphi x^0$ and $\varphi x^j > \varphi x^k$.
A set $G \subseteq \Lattice$ is a \textbf{$\varphi$-Gr\"obner basis}
of $\Lattice$ if for all $b\in\Zset^n$ and for every pair $x,y \in
\Feasible_{\Lattice,b}$, there exists a $\varphi$-reduction path from $x$ to $y$
in $\Graph(\Feasible_{\Lattice,b},G)$.

The following lemma is fundamental to the Project-and-Lift algorithm.
Note that the property that $\ker(\pi_{\{i\}})\cap\Lattice = \{0\}$
for some $i \in \{1,\dots,n\}$ means that each vector in $\Lattice^{\{i\}}$
lifts to a unique vector in $\Lattice$, and thus,
the inverse map $\pi^{-1}_{\{i\}}:\Lattice^{\{i\}} \mapsto \Lattice$ is
well-defined.
Moreover, by linear algebra, there must exist a vector $\omega^i \in \Qset^{n-1}$ 
such that for all $u \in \Lattice^{\{i\}}$, we have
$\omega^i u  = (\pi^{-1}_{\{i\}}(u))_i$.
We will always write such a vector as $\omega^i$, and also, we define
$\bar{\omega}^i=-\omega^i$.

\begin{lem}\label{Lemma:project-and-lift}
Let $i \in \{1,\dots,n\}$ where $\ker(\pi_{\{i\}})\cap\Lattice = \{0\}$,
and let $S \subseteq \Lattice^{\{i\}}$.
The set $S$ is a $\bar{\omega}^i$-Gr\"obner basis of $\Lattice^{\{i\}}$ if and
only if $\pi^{-1}_{\{i\}}(S)$ is a $\bar{e}^i$-Gr\"obner basis of $\Lattice$.
\end{lem}
\boproof
Assume $S$ is a $\bar{\omega}^i$-Gr\"obner basis of $\Lattice^{\{i\}}$.
Let $x,y \in \Feasible_{\Lattice, b}$ for some $b \in \Zset^n$.
We need to show that there is a $\bar{e}^i$-reduction path from $x$ to $y$ in
$\Graph(\Feasible_{\Lattice,b},\pi^{-1}_{\{i\}}(S))$.
Let $\bar{x}=\pi_{\{i\}}(x)$, $\bar{y}=\pi_{\{i\}}(y)$, and
$\bar{b}=\pi_{\{i\}}(b)$.
By assumption, there exists a $\bar{\omega}^i$-reduction path
$(\bar{x}=\bar{x}^0,\dots,\bar{x}^k=\bar{y})$
in $\Graph(\Feasible_{\Lattice^{\{i\}},\bar{b}},S)$.
So, we have $\omega^i\bar{x}^j \ge \omega^i\bar{x}$ or
$\omega^i\bar{x}^j \ge \omega^i\bar{y}$ for all $j$.
We now lift this $\bar{\omega}^i$-reduction path in
$\Graph(\Feasible_{\Lattice^{\{i\}},\bar{b}},S)$ to a $\bar{e}^i$-reduction
path in $\Graph(\Feasible_{\Lattice,b},\pi^{-1}_{\{i\}}(S))$.

Let
$x^j=x+\pi^{-1}_{\{i\}}(\bar{x}^j-\bar{x})
    =y+\pi^{-1}_{\{i\}}(\bar{x}^j-\bar{y})$
for all $j=0,\dots,k$.
Hence, $\pi_{\{i\}}(x^j)=\bar{x}^j$ and
$x^j_i=x_i+\omega^i\bar{x}^j-\omega^i\bar{x}
=y_i+\omega^i\bar{x}^j-\omega^i\bar{y}$,
and so, $x^j_i \ge x_i$ or $x^j_i \ge y_i$.
Also, $x^j-x^{j-1} = \pi^{-1}_{\{i\}}(\bar{x}^j-\bar{x}^{j-1}) \in
\pi^{-1}_{\{i\}}(S)$ for all $j=1,\dots,k$.
Therefore, $(x=x^0,\dots,y^k=y)$ is an
$\bar{e}^i$-reduction path in
$\Graph(\Feasible_{\Lattice,b},\pi^{-1}_{\{i\}}(S))$ as required.

Assume $\pi^{-1}_{\{i\}}(S)$ is a $\bar{e}^i$-Gr\"obner basis of $\Lattice$.
Let $x, y \in \Feasible_{\Lattice^{\{i\}}, b}$ for some $b \in \Zset^{n-1}$, and
let $\gamma = \omega^i(x-y)$. If $\gamma > 0$, then let
$\bar{x} = (x,\gamma)$ and $\bar{y} = (y,0)$,
else let $\bar{x} = (x,0)$ and $\bar{y} = (y,-\gamma)$;
hence, $\bar{x},\bar{y} \in \Feasible_{\Lattice, \bar{b}}$ for some
$\bar{b} \in \Zset^n$,
and $\pi_{\{i\}}(\bar{b}) = b$, $\pi_{\{i\}}(\bar{x}) = x$,
$\pi_{\{i\}}(\bar{y}) = y$, and
$\min\{\bar{x}_i,\bar{y}_i\} = 0$.
By assumption, there exists a $\bar{e}^i$-reduction path
$(\bar{x}={\bar{x}}^0,\dots,{\bar{x}}^k=\bar{y})$
in $\Graph(\Feasible_{\Lattice,\bar{b}},\pi^{-1}_{\{i\}}(S))$.
Let $x^j = \pi_{\{i\}}({\bar{x}}^j)$. So,
$(x=x^0,\dots,x^k=y)$ is a path in $\Graph(\Feasible_{\Lattice^{\{i\}},b},S)$.
Moreover, since $\bar{x}^j_i \ge \bar{x}_i$ or $\bar{x}^j_i \ge \bar{y}_i$
for all $j$, we have $\omega^i x^j \ge \omega^i x$ or 
$\omega^i x^j \ge \omega^i y$ for all $j$.
Therefore, the path is a $\bar{\omega}^i$-reduction path.
\eoproof

By definition, a $\bar{e}^i$-Gr\"obner basis of $\Lattice$ is also
a generating set of $\Lattice$.
On the other hand, a generating set of $\Lattice$ is also a
$\bar{e}^i$-Gr\"obner basis of $\Lattice$.
This follows since, given a generating set of $\Lattice$,
for any $x,y \in \Feasible_{\Lattice,b}$ for any $b$,
there must exist a path from $x-\gamma$ to $y-\gamma$
where $\gamma = x \wedge_{\{i\}} y$, and by translating such a path by $\gamma$,
we get a $\bar{e}^i$-reduction path from $x$ to $y$.
This can also be shown using Lemma \ref{Lemma: saturation}.
So, we arrive at the following corollary.

\begin{cor}
\label{Corollary: lifting a generating set}
Let $i \in \{1,\dots,n\}$ where $\ker(\pi_{\{i\}})\cap\Lattice = \{0\}$,
and let $S \subseteq \Lattice^{\{i\}}$.
The set $S$ is a $\bar{\omega}^i$-Gr\"obner basis of $\Lattice^{\{i\}}$
if and only if $\pi^{-1}_{\{i\}}(S)$ is a generating set of $\Lattice$.
\end{cor}

Given a vector $\varphi \in \Qset^n$,
any $\prec_\varphi$-reduction
path is also a $\varphi$-reduction path, and so, any $\prec_\varphi$-Gr\"obner
basis is also a $\varphi$-Gr\"obner basis.
So, given a set $S \subseteq \Lattice^{\{i\}}$ that generates
$\Lattice^{\{i\}}$, we can compute a $\bar{\omega}^i$-Gr\"obner basis
$S' \subseteq \Lattice^{\{i\}}$ of $\Lattice^{\{i\}}$
by running the completion procedure with respect to $\prec_{\bar{\omega}^i}$ on
$S$.  That is, $S'=\CP(\prec_{\bar{\omega}^i},S)$.
Hence, by Lemma \ref{Lemma:project-and-lift},
$\pi^{-1}_{\{i\}}(S')$ is a generating set of $\Lattice$.

We can apply the above reasoning to compute a generating set of
$\Lattice^{\sigma\setminus\{i\}}$ from a generating set of $\Lattice^\sigma$ for
some $\sigma \subseteq \{1,\dots,n\}$ and $i \in \sigma$.
First, analogously to $\pi_{\{i\}}$ and $\prec_{\bar{\omega}^i}$ in the context
of $\Lattice^{\{i\}}$ and $\Lattice$, we define
$\pi^\sigma_{\{i\}}$ and $\prec^\sigma_{\bar{\omega}^i}$ in the same way except
in the context of $\Lattice^\sigma$ and $\Lattice^{\sigma\setminus\{i\}}$
respectively.

We can now present our Project-and-Lift algorithm.

\begin{alg}{Project-and-Lift algorithm}
\label{Project-and-Lift algorithm}

\underline{Input:} a set $S \subseteq \Lattice$ that spans $\Lattice$.

\vspace{-0.35cm}

\underline{Output:} a generating $G$ set of $\Lattice$

\vspace{0.1cm}

Find a set $\sigma\subseteq\{1,\ldots,n\}$ such that
$\ker(\pi_\sigma) \cap \Lattice = \{0\}$
and $\Lattice^\sigma \cap \Nset^{|\sigmabar|} = \{0\}$.

\vspace{-0.35cm}

Compute a set $G \subseteq \Lattice^\sigma$ such that $G$ is a generating set
of $\Lattice^\sigma$ using $S$.

\vspace{-0.35cm}

\underline{while} $\sigma\neq\emptyset$ \underline{do}

\vspace{-0.35cm}

\hspace{1.0cm}Select $i \in \sigma$

\vspace{-0.35cm}

\hspace{1.0cm}$G:= (\pi^\sigma_{\{i\}})^{-1}(
\CP(\prec^\sigma_{\bar{\omega}^i},G))$

\vspace{-0.35cm}

\hspace{1.0cm}$\sigma:=\sigma\setminus\{i\}$

\vspace{-0.35cm}

\underline{return} $G$.
\end{alg}
\begin{lem}
\label{Lemma: Project-and-Lift algorithm terminates and is correct}
Algorithm \ref{Project-and-Lift algorithm} terminates and satisfies its
specifications.
\end{lem}
\boproof
Algorithm \ref{Project-and-Lift algorithm} terminates, since
Algorithm \ref{Algorithm: Completion procedure} always terminates.

We claim that for each iteration of the algorithm, $G$ is a generating set of
$\Lattice^\sigma$, $\ker(\pi_\sigma) \cap \Lattice = \{0\}$,
and $\Lattice^\sigma \cap \Nset^{|\sigmabar|} = \{0\}$;
therefore, at termination, $G$ is a generating set of $\Lattice$.
This is true for the first iteration, so we assume it is true for
the current iteration.

If $\sigma = \emptyset$, then there is nothing left to do, so assume otherwise.
Since by assumption, $\Lattice^\sigma \cap \Nset^{|\sigmabar|} = \{0\}$ and
$\ker(\pi_\sigma) \cap \Lattice^\sigma = \{0\}$, we must have
$\ker(\pi^\sigma_{\{i\}}) \cap \Lattice^{\sigma\setminus\{i\}} = \{0\}$,
and so, the inverse map
$(\pi^\sigma_{\{i\}})^{-1}: \Lattice^\sigma \rightarrow
\Lattice^{\sigma\setminus\{i\}}$ is well-defined.
Let $i \in \sigma$,
$G' := (\pi^\sigma_{\{i\}})^{-1}(\CP(\prec^\sigma_{\{i\}},G))$, and
$\sigma' := \sigma\setminus\{i\}$.
Then, by Corollary \ref{Corollary: lifting a generating set}, $G'$ is a
generating set of $\Lattice^{\sigma'}$. Also, since $\sigma' \subseteq \sigma$,
we must have $\ker(\pi_{\sigma'}) \cap \Lattice = \{0\}$
and $\Lattice^{\sigma'} \cap \Nset^{|\sigmabar'|} = \{0\}$. Thus, the claim is
true for the next iteration.
\eoproof

In our Project-and-Lift algorithm, we need to find a set
$\sigma\subseteq\{1,\ldots,n\}$
such that $\ker(\pi_\sigma) \cap \Lattice^\sigma = \{0\}$
and $\Lattice^\sigma \cap \Nset^{|\sigmabar|} = \{0\}$, and then, we need to
compute a generating set of $\Lattice^\sigma$.

For our purposes, the larger $\sigma$ the better. However, in general,
finding the largest $\sigma$ is difficult; thus, we use the
following method for finding a \emph{good} $\sigma$. Let $B$ be a
basis for the lattice $\Lattice$ ($\Lattice$ is spanned by the rows of $B$).
Let $k := \rank(B)$. Any $k$ linearly independent columns of $B$ then suffice to
give a set $\sigmabar$ such that every vector in $\Lattice^\sigma$
lifts to a unique vector in $\Lattice$; that is,
$\ker(\pi_\sigma) \cap \Lattice^\sigma = \{0\}$.
Such a set $\sigmabar$ can be found via Gaussian elimination.
If $\Lattice^\sigma\cap\Nset^{|\sigmabar|}\ne\{0\}$, then remove some
$i\in\sigma$ from $\sigma$ ($\sigma:=\sigma\setminus\{i\}$) and recompute
$\Lattice^\sigma\cap\Nset^{|\sigmabar|}$. Continue to do so until
$\Lattice^\sigma\cap\Nset^{|\sigmabar|}=\{0\}$. This procedure must
terminate since $\Lattice\cap\Nset^n=\{0\}$ by assumption.
To check if $\Lattice^\sigma\cap\Nset^{|\sigmabar|}=\{0\}$, we can either solve a
linear programming problem or compute the extreme rays of
$\Lattice^\sigma\cap\Nset^{|\sigmabar|}$ (see for
example \citet{Avis+Fukuda:96,Fukuda+Prodon:96,Hemmecke:Hilbert+Rays}).
In practice, we compute extreme rays using the algorithm in
\citet{Hemmecke:Hilbert+Rays}.

Once we have found such a $\sigma$, we can compute
a generating set of $\Lattice^\sigma$ using either the
Saturation algorithm, the Min-Max algorithm, or any other such algorithm.
In practice and for this paper, we use the Saturation algorithm.
It is possible that there does not exist such a $\sigma$ except the trivial case
where $\sigma = \emptyset$, and so, the Project-and-Lift
algorithm reduces to just the initial phase of computing a generating set of
$\Lattice$ using some other algorithm.
Though in practice, we usually found a non-trivial $\sigma$.
We refer the reader to Section $\ref{Section: What if}$ for a
description of a complete Project-and-Lift algorithm whereby we do not need
another algorithm to start with. 

\begin{exmp}
Consider again the set
$S :=$ \emph{\{(1,-1,-1,-3,-1,2),(1,0,2,-2,-2,1)\}}.
Let $\Lattice$ be the lattice
spanned by $S$. Let $\sigma =$\emph{\{3,4,5,6\}}.
Then, $\ker(\pi_\sigma) \cap \Lattice^\sigma = \{0\}$. 
Note that $\pi_\sigma(S)=$\emph{\{(1,-1),(1,0)\}}.
However, $\Lattice^\sigma\cap\Nset^{|\sigmabar|}\ne\{0\}$.
So, set $\sigma =$ \emph{\{3,4,6\}}.
Now $\pi_\sigma(S)=$\emph{\{(1,-1,-1),(1,0,-2)\}}, and 
$\Lattice^\sigma\cap\Nset^{|\sigmabar|}=\{0\}$.

The set $G = $\emph{\{(0,-1,1),(-1,2,0)\}} is a generating set of
$\Lattice^\sigma$.
We can compute this using the saturation algorithm.
The following table gives the values of $\sigma$, $i$, $\omega^i$,
and $G$ at each stage of the Project-and-Lift algorithm.
\emph{\small{
\begin{center}
\begin{tabular}{|l|l|l|l|l|}
\hline
$\sigma$ & $i$ & $\omega^i$     & $\CP(\prec^\sigma_{\bar{\omega}^i,G})$
& $G:= (\pi^\sigma_{\{i\}})^{-1}(\CP(\prec^\sigma_{\omega^i},G))$ \\
\hline
\hline
\{3,4,6\}& 3 & (2,3,0)      & \{(0,-1,1),(1,-2,0)\}
                            & \{(0,-1,-3,1),(1,-2,-4,0)\} \\
\hline
\{4,6\}  & 4 & (-2,1,0,0)   & \{(0,-1,-3,1),(1,-2,-4,0)\}
                            & \{(0,-1,-3,-1,1),(1,-2,-4,-4,0)\} \\
\hline
\{6\}    & 6 & (1,-1,0,0,0) & \{(0,1,3,1,-1), (-1,1,1,3,1),
                            & \{(0,1,3,1,-1,-1), (-1,1,1,3,1,-2), \\
  &             &           &   (-1,0,-2,2,2), (-1,-1,-5,1,3),
                            &   (-1,0,-2,2,2,-1), (-1,-1,-5,1,3,0), \\
  &             &           &   ( 1,2,8,0,-4)\}
                            &   ( 1,2,8,0,-4,-1)\} \\
\hline
\end{tabular}
\end{center}
}}
The final $G$ is a generating set of $\Lattice$.
Again note that it is not minimal;
the vector \emph{(1,2,8,0,-4,-1)} is not needed and can be removed from $G$,
and $G$ will still be a generating set of $\Lattice$.
\end{exmp}

The concepts of $\sigma$-generating sets of $\Lattice$ and generating sets
of $\Lattice^\sigma$ are, in fact, equivalent.
So, as discussed before, unlike that Saturation
algorithm, the Project-and-Lift algorithm computes $\sigma$-generating sets and
thus does less work than the Saturation algorithm.
\begin{lem}\label{Lemma: gen sets}
Let $\sigma \subseteq \{1,\dots,n\}$ where
$\ker(\pi_\sigma) \cap \Lattice = \{0\}$ and $S \subseteq \Lattice^\sigma$.
The set $S$ is a generating set of $\Lattice^\sigma$ if and only if
$\pi^{-1}_\sigma(S)$ is a $\sigma$-generating set of $\Lattice$.
\end{lem}
\boproof
Recall that a $\sigma$-generating set of $\Lattice$ is a set where for all
for all $b\in\Zset^n$ and for all $x,y \in \Feasible^\sigma_{\Lattice,b}$, there
exists a path from $x$ to $y$ in $\Graph(\Feasible^\sigma_{\Lattice,b},S)$.
Observe that
$\pi_\sigma(\Feasible^\sigma_{\Lattice,b})=
\Feasible_{\Lattice^\sigma,\pi_\sigma(b)}$, and
moreover, a path
in $\Graph(\Feasible^\sigma_{\Lattice,b},S)$ projects to a path in
$\Graph(\Feasible_{\Lattice^\sigma,\pi_\sigma(b)},S)$.
Hence, a $\sigma$-generating set
of $\Lattice$ projects to a generating set of $\Lattice^\sigma$.
So, if $S$ is a $\sigma$-generating set of $\Lattice$, then $\pi_\sigma(S)$
is a generating set of $\Lattice^\sigma$.
Also, assuming $\ker(\pi_\sigma) \cap \Lattice = \{0\}$, if $S$ is a generating
set of $\Lattice^\sigma$, then $\pi^{-1}_\sigma(S)$ is a $\sigma$-generating
set of $\Lattice$.
This follows since a path in $\Graph(\Feasible_{\Lattice^\sigma,b},S)$ can be
lifted to a path in $\Graph(\Feasible^\sigma_{\Lattice,\bar{b}},S)$
where $\pi_\sigma(\bar{b})=b$.
\eoproof

Observe that if $\ker(\pi_\sigma) \cap \Lattice \ne \{0\}$, then a path in 
$\Graph(\Feasible_{\Lattice^\sigma,b},\pi^{-1}_\sigma(S))$ cannot necessarily be
lifted to a path in
$\Graph(\Feasible^\sigma_{\Lattice,\bar{b}},\pi^{-1}_\sigma(S))$
-- the path may become disconnected -- although, we can easily rectify this by
adding a spanning set of the lattice $\ker(\pi_\sigma)\cap \Lattice$ to
$\pi^{-1}_\sigma(S)$.

The Project-and-Lift algorithm has some interesting properties.
As we saw in Lemma \ref{Lemma:project-and-lift},
$\bar{\omega}^i$-reduction paths lift to
$\bar{e}^i$-reduction paths and $\bar{e}^i$-reduction paths project to
$\bar{\omega}^i$-reduction paths. The same holds true for 
$\prec_{\bar{\omega}^i}$-reduction paths and $\prec_{\bar{e}^i}$-reduction
paths, shown in exactly the same way, giving the following lemma.
\begin{lem}\label{Lemma: gb equiv}
Let $i \in \{1,\dots,n\}$ where $\ker(\pi_{\{i\}})\cap\Lattice = \{0\}$,
and let $S \subseteq \Lattice^{\{i\}}$. Let $\prec$ be a term order.
The set $\pi^{-1}_{\{i\}}(S)$ is a $\prec_{\bar{e}^i}$-Gr\"obner basis of
$\Lattice$ if and only if $S$ is a $\prec_{\bar{\omega}^i}$-Gr\"obner basis of
$\Lattice^{\{i\}}$.
\end{lem}
So, during the Project-and-Lift algorithm, we compute a
$\prec_{\bar{\omega}^i}$-Gr\"obner basis for some $i$
and then lift it to a $\prec_{\bar{e}^i}$-Gr\"obner basis. We then compute a
$\prec_{\bar{\omega}^j}$-Gr\"obner basis using some $j\ne i$ and again lift it
to a $\prec_{\bar{e}^j}$Gr\"obner basis, and repeat. So effectively, the
Project-and-Lift algorithm just converts one Gr\"obner basis into another and
lifts to another Gr\"obner basis. We therefore could use a Gr\"obner walk
algorithm to move from one Gr\"obner basis to another
(see \citet{Collart+Kalkbrener+Mall,Fukuda+Jensen+Lauritzen+Thomas}).
We have not yet implemented such an algorithm. It would be interesting to see
its performance.

There are, in fact, two essentially equivalent ways to compute a
$\prec_{\bar{e}^i}$-Gr\"obner basis of $\Lattice$ from a set
$S \subseteq\Lattice^{\{i\}}$ that generates $\Lattice^{\{i\}}$ for some
$i \in \{1,\dots,n\}$ where $\ker(\pi_{\{i\}})\cap\Lattice = \{0\}$,
as needed by the Project-and-Lift algorithm.
Firstly, as we have already seen,
the set $T = \pi^{-1}_{\{i\}}(\CP(\prec_{\bar{w}^i},S))$ is a
$\prec_{\bar{e}^i}$-Gr\"obner basis of $\Lattice$,
but also, the set $T'= \CP(\prec_{\bar{e}^i}, \pi^{-1}_{\{i\}}(S))$ is also a
$\prec_{\bar{e}^i}$-Gr\"obner basis of $\Lattice$.
Essentially, to compute a $\prec_{\bar{e}^i}$-Gr\"obner basis of $\Lattice$,
we do not need a generating set of $\Lattice$, but instead, we only need a 
$\{i\}$-generating set of $\Lattice$.
This follows from the following Lemmas, which are analogous to
Lemmas \ref{Lemma:project-and-lift} and \ref{Lemma: gb equiv} respectively.
\begin{lem}
Let $i \in \{1,\dots,n\}$ where $\ker(\pi_{\{i\}})\cap\Lattice = \{0\}$,
and let $S,T \subseteq \Lattice$ where $S$ is a $\{i\}$-generating set of
$\Lattice$.
The set $T$ is a $\bar{e}^i$-Gr\"obner basis of $\Lattice$ if and only if
for all $b\in\Zset^n$ and for all $x,y \in \Feasible_{\Lattice,b}$ where $x$ and
$y$ are connected in $\Graph(\Feasible_{\Lattice,b},S)$, there exists an
$\bar{e}^i$-reduction path from $x$ to $y$ in
$\Graph(\Feasible_{\Lattice,b},T)$.
\end{lem}
\boproof
The forwards direction must hold by definition.
Conversely, let $x,y \in \Feasible_{\Lattice,b}$ for some $b \in \Zset^n$. Since $S$ is a
$\{i\}$-generating set of $\Lattice$, there must exist $\gamma \in \Nset^n$ where
$\supp(\gamma) \subseteq \{i\}$, such that $x+\gamma$ is connected to $y+\gamma$
in $\Graph(\Feasible_{\Lattice,b+\gamma},S)$. So, by assumption, there exists an
$\bar{e}^i$-reduction path from $x+\gamma$ to $y+\gamma$ in
$\Graph(\Feasible_{\Lattice,b+\gamma},T)$, which translates to a
path from $x$ to $y$ in $\Graph(\Feasible_{\Lattice,b},T)$ since
$\supp(\gamma) \subseteq \{i\}$.
\eoproof

An analogous results holds for $\prec_{\bar{e}^i}$-Gr\"obner bases for similar
reasons.
\begin{lem}
Let $i \in \{1,\dots,n\}$ where $\ker(\pi_{\{i\}})\cap\Lattice = \{0\}$,
and let $S,T \subseteq \Lattice$ where $S$ is a $\{i\}$-generating set of
$\Lattice$.
The set $T$ is a $\prec_{\bar{e}^i}$-Gr\"obner basis of $\Lattice$ if and only if
for all $b\in\Zset^n$ and for all $x,y \in \Feasible_{\Lattice,b}$ where $x$ and
$y$ are connected in $\Graph(\Feasible_{\Lattice,b},S)$, there exists an
$\prec_{\bar{e}^i}$-reduction path from $x$ to $y$ in
$\Graph(\Feasible_{\Lattice,b},T)$.
\end{lem}
Therefore, if $S \subseteq\Lattice^{\{i\}}$ generates $\Lattice^{\{i\}}$,
then $\pi^{-1}_{\{i\}}(S)$ is a $\{i\}$-generating set of $\Lattice$ by Lemma
\ref{Lemma: gen sets}.
So, the set $T'= \CP(\prec_{\bar{e}^i}, \pi^{-1}_{\{i\}}(S))$ is a
$\prec_{\bar{e}^i}$-Gr\"obner basis of $\Lattice$.
Moreover, when computing $\CP(\prec_{\bar{w}^i},S)$ and computing
$\CP(\prec_{\bar{e}^i}, \pi^{-1}_{\{i\}}(S))$, the completion
procedure performs essentially the same sequence of steps producing 
essentially the same output data and intermediate data with the exception that
they perform the computation in different spaces.
These two approaches are thus algorithmically equivalent.

In one iteration, the Saturation algorithm computes
$\CP(\prec_{\bar{e}^i},T)$, in the space $\Lattice$,
for some set $T \subseteq \Lattice$
that is $\sigma$-saturated on some spanning set for some
$\sigma \subseteq \{1,\dots,n\}$ and $i \in \sigmabar$.
On the other hand, in one iteration, the Project-and-Lift effectively computes,
in the space $\Lattice^{\sigma\backslash\{i\}}$, $\CP(\prec_{\bar{e}^i}, T)$
for some set $T \subseteq \Lattice^{\sigma\backslash\{i\}}$ that is a
$\{i\}$-generating set of $\Lattice^\sigma$
for some
$\sigma \subseteq \{1,\dots,n\}$ and $i \in \sigma$.
So, the algorithms are very similar, but the Project-and-Lift algorithm performs
intermediate steps in subspaces whereas the Saturation algorithm performs
intermediate steps in the original space.

\subsection{The ``Lift-and-Project'' algorithm}

The idea behind this algorithm is to lift a spanning set $S$ of
$\Lattice \subseteq \Zset^n$ to a spanning set $S'\subseteq \Zset^{n+1}$ of
$\Lattice'\subseteq \Zset^{n+1}$ in such a way that we can compute a set
$G' \subseteq \Lattice'$ that generates $\Lattice'$ in only one saturation step.
Then, we project $G'$ to $G \subseteq \Lattice$, so that $G$ is a generating set
of $\Lattice$.

Let $S$ be a spanning set of $\Lattice \subseteq \Zset^n$.
Let $S' := \{(u,0): u \in S\} \cup \{(1,\dots,1,-1)\}$, and let
$\Lattice' \subseteq \Zset^{n+1}$ be the lattice spanned by $S'$.
Since the vector $(1,\dots,1,-1)$ is in $S'$,
it follows from Lemma \ref{Lemma: support saturation}, that
if a set $G' \subseteq \Lattice'$ is $\{n+1\}$-saturated on $S'$, then $G'$
is $\{1,\dots,n+1\}$-saturated on $S$, and hence, $G'$ is a
generating set of $\Lattice'$.
Also, since $\Lattice \cap \Nset^n = \{0\}$,
then $\Lattice' \cap \Zset^{n+1}_+ = \{0\}$.
Now, using exactly the same idea behind the Saturation algorithm,
if we let $G' := \CP(\prec_{\bar{e}^{n+1}},S')$, then $G'$ must be a generating
set for the lattice $\Lattice'$ by Lemma \ref{Lemma: saturation}.

So, at the moment, we have a generating set $G'$ for $\Lattice'$, and from this,
we need to extract a generating set of $\Lattice$.
We define the linear map $\rho: \Zset^{n+1} \mapsto \Zset^n$ where
\[\rho(u') := (u'_1+u'_{n+1},u'_2+u'_{n+1},\dots,u'_n+u'_{n+1}).\]
Observe that $\rho$ maps $\Zset^{n+1}$ onto $\Zset^n$, maps $\Lattice'$ onto
$\Lattice$, and maps $\Feasible_{\Lattice',b'}$ onto $\Feasible_{\Lattice,b}$
where $b = \rho(b')$.
Let $G:= \{\rho(u'):u' \in G'\} \setminus \{0\}$.  So, $G \subset \Lattice$,
and we now show that in fact $G$ generates $\Lattice$.
Let $({x'}^0,\dots,{x'}^k)$ be a path in $\Graph(\Feasible_{\Lattice',b'},G')$.
Then, $(\rho({x'}^0),\dots,\rho({x'}^k))$ is a walk from
$\rho({x'}^0)$ to $\rho({x'}^k)$ in $\Graph(\Feasible_{\Lattice,\rho(b)},G)$,
so after removing cycles, we have a path from $\rho({x'}^0)$ to $\rho({x'}^k)$.
Cycles may exist because the kernel of $\rho$ is non-trivial --
$\ker(\rho) = \{(\gamma,\dots,\gamma,-\gamma): \gamma \in \Zset\}$.
Let $x,y \in \Feasible_{\Lattice,b}$ for some $b \in \Zset^n$, and
let $x' := (x,0)$, $y' := (y,0)$, and $b' := (b,0)$; hence,
$\rho(x')=x$, $\rho(y')=y$, and $\rho(b') = b$.
Then, since $G'$ is a generating set of
$\Lattice'$ there must exist a path from $x'$ to $y'$ in
$\Graph(\Feasible_{\Lattice',b'},G')$, and therefore,
there exists a path from $x$ to $y$ in $\Graph(\Feasible_{\Lattice,b},G)$.
Hence, $G$ is a generating set of $\Lattice$. We thus arrive at the
Lift-and-Project algorithm.

\begin{alg}{Lift-and-Project algorithm}
\label{Lift-and-Project Algorithm}

\underline{Input:} a set $S \subseteq \Lattice$ that spans $\Lattice$.

\vspace{-0.35cm}

\underline{Output:} a generating set $G$ of $\Lattice$

\vspace{-0.35cm}

\vspace{0.2cm} $S':=\{(u,0): u \in S\} \cup \{(1,\dots,1,-1)\}$

\vspace{-0.35cm}

$G':= \CP(\prec_{n+1},S')$

\vspace{-0.35cm}

$G := \{\rho(u'):u' \in G'\} \setminus \{0\}$

\vspace{-0.35cm}

\underline{return} $G$.
\end{alg}

To make the algorithm more efficient, we can use a
different additional vector to $(1,\dots,1,-1)$.
By Lemma \ref{Lemma: upper bound on saturations}, we know that given a spanning
set $S$, there exists a $\sigma$ where $|\sigma| \le \lfloor\frac{n}{2}\rfloor$
such that if $T$ is $\sigma$-saturated on $S$, then $T$ is
a generating set of $\Lattice$.
Then, instead of $(1,\dots,1,-1)$, it suffices to use the additional vector
$s_\sigma = \sum_{i\in\sigma} e^i-e^{n+1}$, which has the important property
that $\supp(s_\sigma^+)=\sigma$ and $\supp(s_\sigma^-) = \{n+1\}$.
Recall that $e^i$ is the $i$th unit vector.
Set $S' := \{(u,0): u \in S\} \cup \{s_\sigma\}$, and 
let $\Lattice'$ be the lattice spanned by $S'$.
Then, from Lemma \ref{Lemma: support saturation}, since $s_\sigma \in S'$, if a
set $G' \subseteq \Lattice'$ is $\{n+1\}$-saturated on $S'$, then $G'$ is
$(\sigma\cup\{n+1\})$-saturated on $S'$.
Also, since $\{(u,0): u \in S\} \subseteq S'$, from the proof of Lemma
\ref{Lemma: upper bound on saturations}, it follows that if $G'$ is
$\sigma$-saturated on $S'$, then $G'$ is $\{1,\dots,n\}$-saturated.
Hence, $G'$ is $\{1,\dots,n+1\}$-saturated, and therefore, a generating set of
$\Lattice'$.
So again, we can compute a generating set $G'$ of $\Lattice'$ in
one saturation step.
Also, we similarly define the linear map $\rho_\sigma: \Zset^{n+1} \mapsto \Zset^n$ where
$\rho_\sigma(x') := (x'_1,x'_2,\dots,x'_n)+(\sum_{i\in\sigma}e_i)x_{n+1}$.
Then, $G := \{\rho_\sigma(x'): x' \in G'\}$ is a generating set of $\Lattice$.
As a general rule, the smaller the size of $\sigma$, the faster the algorithm.

\section{What if $\Lattice \cap \Nset^n \ne \{0\}$?}
\label{Section: What if}
If $\Lattice \cap \Nset^n \ne \{0\}$, then computing a generating set of the
lattice $\Lattice$ is actually more straight-forward than otherwise.
The vectors in $\Lattice \cap \Nset^n$ are very useful when
constructing generating sets.

We say that component $i\in\{1,\dots,n\}$ is \textbf{unbounded} if there exists
a $u \in \Lattice \cap \Nset^n$ where $i \in \supp(u)$ and \textbf{bounded}
otherwise.
>From Farkas' lemma, $i$ is unbounded if and only if the linear program
$\max\{x_i: x \equiv 0 \pmod{\Lattice}, x \in \Rset^n_+\}$ is unbounded.
To find a $u \in \Lattice$ such that $u \gneq 0$ and $i \in \supp(u)$,
and so also, to check whether $i$ is unbounded, we can solve a linear program or
compute the extreme rays of $\Lattice \cap \Nset^n$
(see for example \citet{Avis+Fukuda:96,Fukuda+Prodon:96,Hemmecke:Hilbert+Rays}).
Given a term order $\prec$ of $\Lattice$, the order $\prec_{\bar{e}^i}$ is a
term order if and only if $i$ is bounded.

Using the following lemma, we can extend the Saturation algorithm to the more
general case where $\Lattice \cap \Nset^n \ne \{0\}$.
\begin{lem}
\label{Lemma: sat non-negative vector}
Let $S \subseteq \Lattice$. If there exists
$u \in S$ where $ u \in \Lattice \cap \Nset^n$ and $u \ne 0$, then
$S$ is $\supp(u)$-saturated (on $S$).
\end{lem}
\boproof
By definition, $S$ is $\emptyset$-saturated (on $S$).
Since $u \ge 0$, we have $\supp(u^-) = \emptyset$, and so
it follows immediately from Lemma \ref{Lemma: support saturation} that
$S$ is $\supp(u)$-saturated (on $S$).
\eoproof

We can now extend the Saturation algorithm.
Let $\tau \subseteq \{1,\dots,n\}$ be the set of unbounded components,
and let $S$ be a spanning set of $\Lattice$.
Then, for each $i \in \tau$, find a $u \in \Lattice$ such that $u \geq 0$ and
$u_i>0$ and add $u$ to $S$. Or equivalently, find a single $u\geq 0$ such that
$\supp(u) = \tau$ and add it to $S$. Now $S$ is $\tau$-saturated (on $S$) by
Lemma \ref{Lemma: sat non-negative vector}.
So, if a set $T$ is $\taubar$-saturated on $S$, then $T$ is
$\{1,\dots,n\}$-saturated on $S$ by Lemma \ref{Lemma: saturation union}, and so,
$T$ is a generating set of $\Lattice$.
So, we iteratively compute $T := \CP(\prec_{\bar{e}^i},T)$ for every
$i \in \taubar$; then, $T$ is $\taubar$-saturated on $S$ as required.

The Project-and-Lift algorithm can also be extended to the more general case
where $\Lattice \cap \Nset^n \ne \{0\}$. First, we need the following lemma.
\begin{lem}
\label{Lemma: gen non-negative vector}
Let $\sigma \subseteq \{1,\dots,n\}$,
$S \subseteq \Lattice$, and $u \in \Lattice \cap \Nset^n$ and $u \ne 0$.
If $S$ is a $\sigma$-generating set of $\Lattice$ and $u \in S$,
then $S$ is a $(\sigma\setminus\supp(u))$-generating set of
$\Lattice$.
\end{lem}
\boproof
Let $x,y \in \Feasible_{\Lattice,b}$ for some $b \in \Zset^n$.
Since $S$ is a $\sigma$-generating set of $\Lattice$, there exists a path from
$x$ to $y$ in $\Graph(\Feasible^\sigma_{\Lattice,b}, S)$.
This path can transformed into a path in
$\Graph(\Feasible^{(\sigma\setminus\supp(u))}_{\Lattice,b}, S)$
by adding $u$ to the start of the path as many times as necessary and
subtracting $u$ from the end of the path the same number of times.
\eoproof

Let $\sigma \subseteq \{1,\dots,n\}$
where $\ker(\pi_\sigma) \cap \Lattice^\sigma = \{0\}$.
Let $S \subseteq \Lattice^\sigma$ where $S$ is a generating set of
$\Lattice^\sigma$. Also, let $i \in \sigma$. We now show how to construct a
generating set of $\Lattice^{\sigma \setminus \{i\}}$, and thus by induction,
a generating set of $\Lattice$.
Firstly, since $S$ is a generating set of $\Lattice^\sigma$,
$\pi_{\{i\}}^{-1}(S)$ is a $\{i\}$-generating of
$\Lattice^{\sigma \setminus \{i\}}$ from Lemma \ref{Lemma: gen sets}.
If $i$ is unbounded for $\Lattice^{\sigma \setminus \{i\}}$, then there exists
a $u \in \Lattice^{\sigma \setminus \{i\}}$ such that $u \ge 0$ and
$u_i > 0$.  Thus, after adding $u$ to $(\pi^\sigma_{\{i\}})^{-1}(S)$, we
then have a generating set of $\Lattice^{\sigma \setminus \{i\}}$.
If $i$ is bounded, then compute
$S :=\CP(\prec_{\bar{\omega}^i},S)$, and
$(\pi^\sigma_{\{i\}})^{-1}(S)$
is then a generating set of $\Lattice^{\sigma \setminus \{i\}}$.

We first need to find an initial $\sigma$ and $S$ such that 
$\ker(\pi_\sigma) \cap \Lattice^\sigma = \{0\}$ and $S$ is a generating set of
$\Lattice^\sigma$.
Let $B$ be a lattice basis of $\Lattice$ where the rows of $B$ span $\Lattice$.
Let $\sigmabar$ be any $\rank(B)$ linearly independent columns of $B$.
Let $S = \pi_\sigma(B)$.
Then, every vector in $\Lattice^\sigma$ lifts to a unique vector in $\Lattice$.
Then, computing a $u \in \Lattice^\sigma$ such that $u > 0$ can be
done by Gaussian elimination. After adding $u$ to $S$,
$S$ is a generating set of $\Lattice^\sigma$ as required.

We can also extend the Lift-and-Project algorithm in a similar way. As above, we
can find a set $S$ such that $S$ is $\tau$-saturated (on $S$).
The set \[S' := \{(u,0): u \in S\} \cup \{\sum_{i\in\taubar} e_i-e_{n+1}\}\]
is also $\tau$-saturated (on $S'$) for the same reasons.
If a set $T' \subseteq \Lattice'$ is $\{n+1\}$-saturated on $S'$, then $T'$ is
$(\taubar\cup\{n+1\})$-saturated on $S'$
by Lemma \ref{Lemma: sat non-negative vector}, and so, 
by Lemma \ref{Lemma: saturation union}, $T'$ is $\{1,...,n+1\}$-saturated on
$S'$ since $S'$ is $\tau$-saturated (on $S'$);
thus, $T'$ is then a generating set of $\Lattice'$.
Hence, in one saturation step, we can compute a generating set of
$\Lattice'$. Note that the component $n+1$ is bounded by construction.
Then, the set $T := \{\rho_{\taubar}(u'):u'\in T'\}$ is a generating set of
$\Lattice$.

\section{Speeding-up the Completion Procedure}
\label{Section: Speeding-up the Completion Procedure}

Finally, before presenting computational experience, we talk
about ways in which the key algorithm,
Algorithm \ref{Algorithm: Completion procedure}, can be improved.
This leads us to the critical pair criteria.

Algorithm \ref{Algorithm: Completion procedure} has to test for a
reduction path between $x^{(u,v)}$ and $y^{(u,v)}$ for all critical pairs
$C := \{(u,v): u,v \in G\}$.
In the case of lattice ideals, computational profiling
shows that this is the most time consuming part of the computation.
So, we wish to reduce the number of critical pairs that we test, and
avoid this expensive test as often as possible.
We present three criteria that can reduce the number of critical pairs
that need to be tested.

Criteria $1$ and $3$ (see \citet{Buchberger:79,Buchberger:85,Gebauer+Moeller:88})
are translated from the theory of Gr\"obner bases into a geometric context.
Criterion $2$ is specific to lattice ideals and corresponds to using the
homogeneous Buchberger algorithm
\citep{Caboara+Kreuzer+Robbiano:03,Traverso:97},
but we give a slightly more general result.
Note that all three criteria can be applied simultaneously.

\subsubsection*{Criterion $1$: The Disjoint-Positive-Support criterion}

For a pair $u,v \in G$, the Disjoint-Positive-Support criterion is a simple and
quick test for a $\succ$-reduction path from $x^{(u,v)}$ to $y^{(u,v)}$.
So, using this quick test for a $\succ$-reduction path,
we can sometimes avoid the more expensive test.

Given $u,v \in G$, if $\supp(u^+) \cap \supp(v^+) = \emptyset$, then
there exists a simple $\succ$-reduction path from $x^{(u,v)}$ to $y^{(u,v)}$
using $u$ and $v$ in reverse order (see Figure \ref{Figure: Criterion 1}).

\input{Figure_Criterion_1}

\subsubsection*{Criterion $2$: The Cancellation criterion}

Let $G$ be a generating set of $\Lattice$. 
If $\supp(x^{(u,v)})\cap\supp(y^{(u,v)}) \ne \emptyset$ for some $u,v \in G$
(or equivalently, $\supp(u^-)\cap\supp(v^-) \ne \emptyset$), then we do not need
to check for a $\succ$-reduction path from $x^{(u,v)}$ to $y^{(u,v)}$ (we can
remove the pair $(u,v)$ from $C$).

To show that this criteria holds, we need the concept of a grading.
Let $w \in \Qset^n$. If $w x=w y$ for all $x,y \in
\Feasible_{\Lattice,b}$ for all $b \in \Zset^n$, then we call $w$ a
grading of $\Lattice$, and we define $\deg_w(\Feasible_{\Lattice,b})
:= w b$ called the $w$-\textbf{degree} of
$\Feasible_{\Lattice,b}$. Importantly, if $\Lattice \cap \Nset^n =
\{0\}$ (which we assume), then it follows from Farkas' lemma that
there exists a strictly positive grading $w \in \Qset^n_+$ of
$\Lattice$.

First, we prove an analogous result to Corollary \ref{Corollary:
Test for reduction paths of critical paths is sufficient}.
\begin{lem}
A set $G\subseteq\Lattice_\succ$ is a $\prec$-Gr\"obner basis of
$\Lattice$ if and only if $G$ is a
generating set of $\Lattice$ and if for every $\succ$-critical path $(x,z,y)$ in
$\Graph(\Feasible_{\Lattice,b},G)$ for all $b\in\Zset^n$ where
$\supp(x)\cap\supp(y)=\emptyset$, there exists a
$\succ$-reduction path between $x$ and $y$ in
$\Graph(\Feasible_{\Lattice,b},G)$.
\end{lem}
\boproof The forwards implication follows from Corollary
\ref{Corollary: Test for reduction paths of critical paths is sufficient}.
For the backwards implication, we need to show that for every
$\succ$-critical path $(x,y,z)$ where $\supp(x)\cap\supp(y) \ne \emptyset$,
there exists a $\succ$-reduction path from $x$ to $y$ in
$\Graph(\Feasible_{\Lattice,x},G)$, in which case, there is a
$\succ$-reduction path for all $\succ$-critical paths, and so by Corollary
\ref{Corollary: Test for reduction paths of critical paths is sufficient},
$G$ is a Gr\"obner basis. Assume on the contrary that this is not
the case. Let $w$ be a strictly positive grading of $\Lattice$.
Among all such $\succ$-critical paths $(x,z,y)$ where
$\supp(x)\cap\supp(y) \ne \emptyset$ and there is no $\succ$-reduction path from
$x$ to $y$, choose a $\succ$-critical path $(x,z,y)$ such that
$\deg_w(\Feasible_{\Lattice,x})$ is minimal.
Let $\gamma := x \wedge y$,
$\bar{x}:=x-\gamma$, and $\bar{y}:=y-\gamma$. Note that
$\gamma \ne 0$ since $\supp(x)\cap\supp(y) \ne \emptyset$. Because $G$ is a
generating set of $\Lattice$, there must exist a path from $\bar{x}$
to $\bar{y}$ in $\Graph(\Feasible_{\Lattice,\bar{x}},G)$. Also,
since $w$ is strictly positive,
$\deg_w(\Feasible_{\Lattice,\bar{x}}) < \deg_w(\Feasible_{\Lattice,x})$;
therefore, by the minimality assumption on $\deg_w(\Feasible_{\Lattice,x})$, we
can now conclude that for all $\succ$-critical paths in
$\Graph(\Feasible_{\Lattice,\bar{x}},G)$
there exists a $\succ$-reduction path. Consequently, by Lemma \ref{Lemma:
Test for reduction paths of critical paths is sufficient for each fiber},
there exists a $\succ$-reduction path between $\bar{x}$ and $\bar{y}$ in
$\Graph(\Feasible_{\Lattice,\bar{x}},G)$. This $\succ$-reduction path,
however, can be translated by $\gamma$ to a $\succ$-reduction path from $x$
to $y$ in $\Graph(\Feasible_{\Lattice,x},G)$ (see Figure
\ref{Figure: Criterion 2}a). But this contradicts our assumption
that there is no such path between $x$ and $y$. \eoproof

\input{Figure_Criterion_2}

Now, for all $u,v \in G$, if $\supp(x^{(u,v)})\cap\supp(y^{(u,v)})\ne\emptyset$,
then $\supp(x) \cap \supp(y) \ne \emptyset$ for all $\succ$-critical paths
$(x,z,y)$ for $(u,v)$.
Using this observation, we arrive at an analogous result to
Corollary \ref{Corollary: Test for reduction paths of minimal
critical paths is sufficient}.

\begin{cor} \label{Corollary: Criterion 2}
Let $G \subseteq \Lattice$ be a generating set of $\Lattice$; then,
$G$ is a $\prec$-Gr\"obner basis of $\Lattice$ if
and only if for each pair $u,v\in G$ where
$\supp(x^{(u,v)})\cap\supp(y^{(u,v)}) = \emptyset$, there exists a
$\succ$-reduction path between $x^{(u,v)}$ and $y^{(u,v)}$ in
$\Graph(\Feasible_{\Lattice,z^{(u,v)}},G)$.
\end{cor}

We can extend these results further leading to a more powerful
elimination criterion.  Let $u,v \in G$. We say the pair $(u,v)$ satisfies
Criterion $2$ if there exists $x',y' \in \Feasible_{\Lattice,z^{(u,v)}}$
such that there exists a $\succ$-decreasing path in
$\Graph(\Feasible_{\Lattice,z^{(u,v)}},G)$ from $x^{(u,v)}$ to $x'$
and from $y^{(u,v)}$ to $y'$, and $\supp(x') \cap \supp(y') \ne
\emptyset$. Importantly, if $(u,v)$ satisfies Criterion 2, then we do
not have to test for a $\succ$-reduction path from $x^{(u,v)}$ to $y^{(u,v)}$.
Thus, we arrive at an extension of Corollary \ref{Corollary:
Criterion 2}. Observe that the previous results are just a special
case where $x'=x^{(u,v)}$ and $y'=y^{(u,v)}$.

\begin{lem} \label{Lemma: Criterion 2}
Let $G \subseteq \Lattice$ be a generating set of $\Lattice$; then,
$G$ is a $\prec$-Gr\"obner basis of $\Lattice$ if
and only if for each pair $u,v\in G$ where $(u,v)$ does not satisfy
Criterion 2, there exists a $\succ$-reduction path between $x^{(u,v)}$ and
$y^{(u,v)}$ in $\Graph(\Feasible_{\Lattice,z^{(u,v)}},G)$.
\end{lem}

If there is a $\succ$-reduction path from $x'$ to $y'$, then there exists a
$\succ$-reduction path from $x^{(u,v)}$ to $y^{(u,v)}$.
Since $\supp(x') \cap \supp(y') \ne \emptyset$, then
$\gamma = x' \wedge y' \ne 0$.
Let $\bar{x} = x'-\gamma$ and $\bar{y} = y'-\gamma$.
So, if there exists a
$\succ$-reduction path from $\bar{x}$ to $\bar{y}$, then there must exist a
$\succ$-reduction path from $x'$ to $y'$, and therefore also, there must
exist a $\succ$-reduction path from $x^{(u,v)}$ to $y^{(u,v)}$ (see Figure
\ref{Figure: Criterion 2}b). Again, we let $w$ be a strictly
positive grading of $\Lattice$, and so similarly to above,
$\deg_w(\Feasible_{\Lattice,\bar{x}}) <
\deg_w(\Feasible_{\Lattice,z^{(u,v)}})$. So, the proof of Lemma
\ref{Lemma: Criterion 2} is essentially as before.

For a pair $u,v \in G$, Criterion 2 can be checked not only before
we search for a $\succ$-reduction path from $x^{(u,v)}$ to $y^{(u,v)}$ but also
while searching for a $\succ$-reduction path. When searching for a
$\succ$-reduction path, we construct a $\succ$-decreasing path from $x^{(u,v)}$
to $\NF(x^{(u,v)},G)$ and a $\succ$-decreasing path from $y^{(u,v)}$ to
$\NF(y^{(u,v)},G)$.
Therefore, we can take any point $x'$ on the $\succ$-decreasing path from
$x^{(u,v)}$ to $\NF(x^{(u,v)},G)$ and any point $y'$ on the decreasing
path from $y^{(u,v)}$ to $\NF(y^{(u,v)},G)$ and check Criterion 2, that
is, we check if $\supp(x') \cap \supp(y') \neq \emptyset$. If this
is true, then we can eliminate $(u,v)$.

We wish to point out explicitly here that Criterion $2$ can be
applied without choosing the vector pairs $u,v \in G$ in a
particular order during Algorithm \ref{Algorithm: Completion procedure}.
In fact, when running Algorithm \ref{Algorithm: Completion procedure},
if we apply Criterion 2 to eliminate a pair
$u,v \in G$, it does not necessarily mean that there is a $\succ$-reduction
path from $x^{(u,v)}$ to $y^{(u,v)}$ in
$\Graph(\Feasible_{\Lattice,z^{(u,v)}},G)$ at that particular
point in time in the algorithm but instead that a $\succ$-reduction path
will exist when the algorithm terminates. This approach is in
contrast to existing approaches that use the homogeneous Buchberger algorithm
to compute a Gr\"obner basis whereby vector pairs $u,v \in G$ must be chosen
in an order compatible with increasing
$\deg_w(\Feasible_{\Lattice,z^{(u,v)}})$ for some strictly positive
grading $w$. This can be computationally costly.
When we use these existing approaches, if a pair $(u,v)$ is
eliminated by Criterion 2, then it is necessarily the case that
there already exists a $\succ$-reduction path from $x^{(u,v)}$ to $y^{(u,v)}$.

Since we need a generating set of $\Lattice$ for Criterion $2$, we cannot apply
Criterion $2$ during the Saturation algorithm (Algorithm \ref{Saturation
algorithm}), and also, we cannot apply Criterion $2$ when
$\Lattice \cap \Nset^n \ne \{0\}$.
However, we can apply a less strict version.
Given $u,v \in G$ and $\tau \subseteq \{1,\dots,n\}$, we say that $(u,v)$
satisfies Criterion $2$ with respect to $\tau$, if there exists $x',y' \in
\Feasible_{\Lattice,z^{(u,v)}}$ such that there exists a decreasing
path in $\Graph(\Feasible_{\Lattice,z^{(u,v)}},G)$ from $x^{(u,v)}$
to $x'$ and from $y^{(u,v)}$ to $y'$, and $\supp(x') \cap \supp(y')
\cap \tau \ne \emptyset$.

Let $S,T \subseteq \Lattice$ where $S$ spans $\Lattice$ and $T$ is a
$\sigma$-saturated set on $S$ for some $\sigma \subseteq \{1,\dots,n\}$.
During an iteration of the Saturation algorithm \ref{Saturation algorithm},
we compute a $(\sigma\cup\{i\})$-saturated set of $S$, by computing
$\CP(\prec_{\bar{e}^i},T)$.
While computing $\CP(\prec_{\bar{e}^i},T)$ here,
we may apply Criterion 2 with respect to $\sigma$. For an
algebraic proof of this, see \citet{Bigatti+Lascala+Robbiano:99}.

Also, we can use Criterion $2$ when $\Lattice \cap \Nset^n \ne \{0\}$ if we have
a generating set of $\Lattice$.  Let $\tau$ be the set of bounded components.
Then, we may apply Criterion $2$ with respect to $\tau$.  Moreover, if we do
not have a generating set and we are running the Saturation algorithm when
$\Lattice \cap \Nset^n \ne \{0\}$, we may apply Criterion 2 with respect to
$\sigma \cap \tau$.

\subsubsection*{Criterion $3$: The $(u,v,w)$ criterion}

Before presenting the $(u,v,w)$ criterion, we need a another result, Lemma
\ref{Lemma: Test for bounded paths of critical paths is sufficient for
each fiber}, that is a less strict version of Lemma
\ref{Lemma: Test for reduction paths of critical paths is sufficient
for each fiber}.
First, we need to define a new type of path. A path
$(x^0,\dots,x^k)$ is $z$-\textbf{bounded} (with respect to $\prec$)
if $x^i \prec z$
for all $i=0,\dots,k$. So, $z$ is a strict upper bound on the path.
Note that for a $\succ$-critical path $(x,z,y)$, a $\succ$-reduction path from
$x$ to $y$ is a $z$-bounded path.
\begin{lem}
\label{Lemma: Test for bounded paths of critical paths is sufficient for
each fiber} Let $b\in\Zset^n$, $x,y\in \Feasible_{\Lattice,b}$, and let
$G\subseteq\Lattice_\succ$ where there is a path between
$x$ and $y$ in $\Graph(\Feasible_{\Lattice,b},G)$. If there
exists a $z'$-bounded path between $x'$ and $y'$ for every $\succ$-critical path
$(x',z',y')$ in $\Graph(\Feasible_{\Lattice,b},G)$, then there
exists a $\succ$-reduction path between $x$ and $y$ in
$\Graph(\Feasible_{\Lattice,b},G)$.
\end{lem}
If we now re-examine the proof of Lemma
\ref{Lemma: Test for bounded paths of critical paths is sufficient for
each fiber}, we find that we only need $z'$-bounded paths between $x'$ and $y'$
for every $\succ$-critical path $(x',z',y')$ in
$\Graph(\Feasible_{\Lattice,b},G)$,
and that, a $\succ$-reduction path from $x'$ and $y'$ is more than we need.
The proof proceeds in the same way as Lemma \ref{Lemma:
Test for reduction paths of critical paths is sufficient for each fiber}.

>From Lemma \ref{Lemma: Test for bounded paths of critical paths is
sufficient for each fiber}, we arrive at an analogous result to
Corollary \ref{Corollary: Test for reduction paths of minimal
critical paths is sufficient}.
\begin{cor}
\label{Corollary: Test for bounded paths of minimal critical paths is
sufficient} A set $G\subseteq\Lattice_\succ$ is a $\prec$-Gr\"obner
basis of $\Lattice$ if and only if $G$ is a
generating set of $\Lattice$ and if for each pair $u,v\in G$, there
exists a $z^{(u,v)}$-bounded path between $x^{(u,v)}$ and $y^{(u,v)}$ in
$\Graph(\Feasible_{\Lattice,z^{(u,v)}},G)$.
\end{cor}

Corollary \ref{Corollary: Test for bounded paths of minimal critical paths
is sufficient} does not fundamentally change Algorithm
\ref{Algorithm: Completion procedure} since to test for a
$z^{(u,v)}$-bounded path from $x^{(u,v)}$ to $y^{(u,v)}$, we still test for a
$\succ$-reduction path from $x^{(u,v)}$ to $y^{(u,v)}$ which is a
$z^{(u,v)}$-bounded path.
However, we can use Corollary
\ref{Corollary: Test for bounded paths of minimal critical paths is
sufficient} to reduce the number of critical pairs $u,v \in G$ for
which we need to compute a $\succ$-reduction path.

Now, we are able to present the $(u,v,w)$ criterion.
Let $u,v,w\in G$ where $z^{(u,v)} \geq w^+$
(or equivalently, $z^{(u,v)} \geq z^{(u,w)}$ and $z^{(u,v)} \geq z^{(w,v)}$),
and let $\bar{z} = z^{(u,v)}-w$.
Then, a $z^{(u,v)}$-bounded path from $x^{(u,w)}$ to $\bar{z}$,
and a $z^{(u,v)}$-bounded path from $\bar{z}$ to $y^{(w,v)}$
combine to form a $z^{(u,v)}$-bounded path from $x^{(u,v)}$ to $y^{(u,v)}$.
Moreover, $(x^{(u,v)},z^{(u,v)},\bar{z})$ is a $\succ$-critical path for $(u,w)$
and $(\bar{z},z^{(u,v)},y^{(u,v)})$ is a $\succ$-critical path for $(w,v)$
(see Figure \ref{Figure: Criterion 3}).
Therefore, a $z^{(u,w)}$-bounded path from $x^{(u,w)}$ to
$y^{(u,w)}$ and a $z^{(w,v)}$-bounded path from $x^{(w,v)}$ to $y^{(w,v)}$
combine to form a $z^{(u,v)}$-bounded path from $x^{(u,v)}$ to $y^{(u,v)}$,
and so, we can remove $(u,v)$ from $C$.

\input{Figure_Criterion_3}

Note that in Figure \ref{Figure: Criterion 3}a,
a $\succ$-reduction path from $x^{(u,v)}$ to $\bar{z}$ and a
$\succ$-reduction path from $\bar{z}$ to $y^{(u,v)}$ do combine to give a
$\succ$-reduction path from $x^{(u,v)}$ to $y^{(u,v)}$;
however, this is not the case in Figure \ref{Figure: Criterion 3}b which is why
we need the concept of bounded paths.

We can extend the previous result. Let $u,v,\in G$, and $w^1,\dots,w^k \in G$
where $z^{(u,v)} \geq {(w^i)}^+$ for all $i=1,\dots,k$. If there exists a
bounded path for the critical pairs $(u,w^1)$, $(w^k,v)$, and $(w^i,w^{i+1})$
for all $i=1,\dots,k-1$, then there is a bounded path for $(u,v)$.
However, note that this can also be implied by a
bounded path for $(u,w^i)$ and $(w^i,v)$ for any $i=1,\dots,k$.

Unfortunately, we cannot just remove from $C$ all pairs $u,v \in G$ where there
exists a $w \in G$ such that $z^{(u,v)} \ge w^+$. It may happen that in
addition to $z^{(u,v)} \geq w^+$, we also have $z^{(u,w)} \ge v^+$, in
which case, we would eliminate both the pairs $(u,v)$ and $(u,w)$
leaving only $(v,w)$ which is not sufficient. Moreover, at the same
time, we may also have $z^{(w,v)} \ge u^+$, and we would eliminate all
three pairs. To avoid these circular relationships, Gebauer and
M\"oller \citet{Gebauer+Moeller:88} devised the following critical pair
elimination criteria which we use in practice in 4ti2 v1.2.

Let $G = \{u^1,u^2,\dots,u^{|G|}\}$, and let $u^i,u^j \in G$ where
$i < j$. We define that the pair $(u^i,u^j)$ satisfies Criterion 3 if
there exists $u^k \in G$ such that one of the following conditions
hold:
\begin{enumerate}
\item $z^{(u^i,u^j)} \gneq z^{(u^i,u^k)}$, and
$z^{(u^i,u^j)} \gneq z^{(u^j,u^k)}$;
\item $z^{(u^i,u^j)} = z^{(u^i,u^k)}$, $z^{(u^i,u^j)} \gneq z^{(u^j,u^k)}$,
and $k<j$;
\item $z^{(u^i,u^j)} \gneq z^{(u^i,u^k)}$, $z^{(u^i,u^j)} = z^{(u^j,u^k)}$,
and $k<i$; or
\item $z^{(u^i,u^j)} = z^{(u^i,u^k)} = z^{(u^j,u^k)}$, and $k < i < j$.
\end{enumerate}
So, if a pair $(u^i,u^j)$ satisfies Criterion 3, we can eliminate
it. For example, if $G = \{u^1,u^2,u^3\}$ where
$z^{u^1u^2}=z^{u^1u^3}\gneq z^{u^2u^3}$, then applying Criterion 3
to all three pairs $(u^1,u^2)$, $(u^1,u^3)$, and $(u^2,u^3)$ would
eliminate only $(u^1,u^3)$.

After eliminating all pairs that satisfy Criterion 3, we are left
with a set of critical pairs
$C' \subseteq C = \{(u,v): u,v \in G\}$ such that if there exists a
$z^{(u',v')}$-bounded path from $x^{(u',v')}$ to $y^{(u',v')}$ for all
$(u',v') \in C'$, then there exists a $z^{(u,v)}$-bounded path from
$x^{(u,v)}$ to $y^{(u,v)}$ for all $(u,v) \in C$.
However, this set of pairs may not be minimal.
In \citet{Caboara+Kreuzer+Robbiano:03}, Caboara, Kreuzer, and
Robbiano describe an algebraic algorithm for computing a minimal set
of critical pairs with computational results. Their computational
results show that the Gebauer and M\"oller criteria give a good
approximation to the minimal set of critical pairs. We found that
the Gebauer and M\"oller criteria were sufficient for our
computations.

\section{The $4\times 4\times 4$-challenge}
\label{Section: The 4x4x4-challenge}

The challenge posed by Seth Sullivant amounts to checking whether a
given set of $145,512$ integer vectors in $\Zset^{64}$ is a Markov
basis for the statistical model of $4\times 4\times 4$ contingency
tables with $2$-marginals. If $x=(x_{ijk})_{i,j,k=1,\ldots,4}$
denotes a $4\times 4\times 4$ array of integer numbers, the defining
equations for the sampling moves are
\begin{eqnarray*}
\sum_{i=1}^4 x_{ijk} = 0 & & \text{ for } j,k=1,\ldots,4,\\
\sum_{j=1}^4 x_{ijk} = 0 & & \text{ for } i,k=1,\ldots,4,\\
\sum_{k=1}^4 x_{ijk} = 0 & & \text{ for } i,j=1,\ldots,4.\\
\end{eqnarray*}
This leads to a problem matrix $A_{444}\in\Zset^{48\times 64}$ of rank
$37$ and $\Lattice_{A_{444}}=\{z:A_{444}z=0,z\in\Zset^{64}\}$. Note that
the $145,512$ vectors in the conjectured Markov basis fall into $14$
equivalence classes under the natural underlying symmetry group
$S_4\times S_4\times S_4\times S_3$.

In \citet{Aoki+Takemura:03}, Aoki and Takemura have computed these
$14$ symmetry classes via an analysis of sign patterns and under
exploitation of symmetry. They claimed that the corresponding
$145,512$ vectors form the unique inclusion-minimal Markov basis of
$A_{444}$.

Using our Project-and-Lift algorithm, however, we have computed the
Markov basis from the problem matrix $A_{444}$ within less than $7$
days on a Sun Fire V890 Ultra Sparc IV processor with 1200 MHz. Note
that the symmetry of the problem was \emph{not} used by the
algorithm. This leaves room for a further significant speed-up. Our
computation produced $148,968$ vectors; that is, there are
additionally $3,456$ Markov basis elements. These vectors form a
single equivalence class under $S_4\times S_4\times S_4\times S_3$
of a norm $28$ vector $z_{15}$ (or equivalently, of a degree $14$
binomial).

A quick check via a Hilbert basis computation with \FourTiTwo\ shows
that these Markov basis elements are indispensable, since
$\{z\in\Zset_+^{64}:A_{444}z=A_{444}z_{15}^+\}=\{z_{15}^+,z_{15}^-\}$.
As also all the other $145,512$ Markov basis elements were
indispensable, the Markov basis of $4\times 4\times 4$ contingency
tables with $2$-marginals is indeed unique. At least this claim can
be saved from \citet{Aoki+Takemura:03}, although we have finally
given a computational proof. Here is the list of the $15$ orbit
representatives, written as binomials:

{\small \noindent
\begin{enumerate}
\item $x_{111}x_{144}x_{414}x_{441}-x_{114}x_{141}x_{411}x_{444}$
\item $x_{111}x_{144}x_{334}x_{341}x_{414}x_{431}-x_{114}x_{141}x_{331}x_{344}x_{411}x_{434}$
\item
$x_{111}x_{122}x_{134}x_{143}x_{414}x_{423}x_{432}x_{441}-
x_{114}x_{123}x_{132}x_{141}x_{411}x_{422}x_{434}x_{443}$
\item
$x_{111}x_{144}x_{324}x_{333}x_{341}x_{414}x_{423}x_{431}-
x_{114}x_{141}x_{323}x_{331}x_{344}x_{411}x_{424}x_{433}$
\item
$x_{111}x_{144}x_{234}x_{243}x_{323}x_{341}x_{414}x_{421}x_{433}-
x_{114}x_{141}x_{233}x_{244}x_{321}x_{343}x_{411}x_{423}x_{434}$
\item
$x_{111}x_{122}x_{133}x_{144}x_{324}x_{332}x_{341}x_{414}x_{423}x_{431}-\\
x_{114}x_{123}x_{132}x_{141}x_{322}x_{331}x_{344}x_{411}x_{424}x_{433}$
\item
$x_{111}x_{144}x_{222}x_{234}x_{243}x_{323}x_{341}x_{414}x_{421}x_{432}-\\
x_{114}x_{141}x_{223}x_{232}x_{244}x_{321}x_{343}x_{411}x_{422}x_{434}$
\item
$x_{111}x_{144}x_{222}x_{233}x_{324}x_{332}x_{341}x_{414}x_{423}x_{431}-\\
x_{114}x_{141}x_{223}x_{232}x_{322}x_{331}x_{344}x_{411}x_{424}x_{433}$
\item
$x_{111}x_{112}x_{133}x_{144}x_{223}x_{224}x_{232}x_{241}x_{314}x_{322}x_{413}x_{421}-\\
x_{113}x_{114}x_{132}x_{141}x_{221}x_{222}x_{233}x_{244}x_{312}x_{324}x_{411}x_{423}$
\item
$x_{111}x_{112}x_{133}x_{144}x_{224}x_{232}x_{243}x_{313}x_{322}x_{341}x_{414}x_{421}-\\
x_{113}x_{114}x_{132}x_{141}x_{222}x_{233}x_{244}x_{312}x_{321}x_{343}x_{411}x_{424}$
\item
$x_{111}x_{134}x_{143}x_{222}x_{233}x_{241}x_{314}x_{323}x_{342}x_{412}x_{424}x_{431}-\\
x_{114}x_{133}x_{141}x_{223}x_{231}x_{242}x_{312}x_{324}x_{343}x_{411}x_{422}x_{434}$
\item
$x_{111}x_{134}x_{143}x_{224}x_{232}x_{241}x_{314}x_{323}x_{342}x_{412}x_{421}x_{433}-\\
x_{114}x_{133}x_{141}x_{221}x_{234}x_{242}x_{312}x_{324}x_{343}x_{411}x_{423}x_{432}$
\item
$x_{111}^2x_{124}x_{133}x_{144}x_{214}x_{223}x_{242}x_{313}x_{332}x_{341}x_{414}x_{424}x_{431}-\\
x_{114}^2x_{123}x_{131}x_{141}x_{213}x_{222}x_{244}x_{311}x_{333}x_{342}x_{411}x_{424}x_{432}$
\item
$x_{111}^2x_{124}x_{133}x_{144}x_{214}x_{232}x_{243}x_{312}x_{323}x_{341}x_{414}x_{422}x_{431}-\\
x_{114}^2x_{123}x_{131}x_{141}x_{212}x_{233}x_{244}x_{311}x_{322}x_{343}x_{411}x_{424}x_{432}$
\item
$x_{111}^2x_{133}x_{144}x_{223}x_{224}x_{232}x_{242}x_{313}x_{322}x_{341}x_{414}x_{422}x_{431}-\\
x_{113}x_{114}x_{131}x_{141}x_{222}^2x_{233}x_{244}x_{311}x_{323}x_{342}x_{411}x_{424}x_{432}$
\end{enumerate}
}

\section{Computational experience}
\label{Section: Computational experience}

We now compare the implementation of our new algorithm in
\FourTiTwo\ v.1.2 \citep{4ti2} with the implementation of the
Saturation algorithm \citep{Hosten+Sturmfels} and the Lift-and-Project
algorithm \citep{Bigatti+Lascala+Robbiano:99} in Singular v3.0.0
\citep{Singular} (algorithmic options `hs' and `blr') and in CoCoA
4.2 \citep{CoCoA} (functions `Toric' and `Toric.Sequential').

\begin{center}
\begin{tabular}{|l|l|l|l|}
\hline
Name & Software & Function & Algorithm\\
\hline Sing-blr & Singular v3.0.0 & toric, option ``blr'' &
Lift-and-Project \\ 
Sing-hs & Singular v3.0.0 & toric, option ``hs'' &
Saturation \\ 
CoCoA-t & CoCoA v4.2 & Toric & Lift-and-Project \\
CoCoA-ts & CoCoA v4.2& Toric.Sequential & Saturation \\
P\&L & \FourTiTwo\ v1.2 & groebner & Project-and-Lift \\
\FourTiTwo\text{-gra} & \FourTiTwo\ v1.2 & graver & Graver basis 
\citep{Hemmecke:SymmGraver}\\
\hline
\end{tabular}
\end{center}

The first $4$ problems correspond to three-way tables with
$2$-marginals, whereas $K4$ and $K5$ correspond to the binary models
on the complete graphs $K_4$ and $K_5$, respectively. The problem
\texttt{s-magic333} is taken from an application in
\citet{Ahmed+DeLoera+Hemmecke:2003} and computes the relations among
the $66$ elements of the Hilbert basis elements of $3\times 3\times
3$ semi-magic hypercubes. The example \texttt{grin} is taken from
\citet{Hosten+Sturmfels}, while the examples \texttt{hppi10-hppi14} correspond
to the computation of homogeneous primitive partition identities,
see for example Chapters $6$ and $7$ in \citet{Sturmfels:96}.
Finally, the examples \texttt{cuww1-cuww5} arise from knapsack problems
presented in \citet{Cornuejols+Urbaniak+Weismantel+Wolsey}.

\vspace{0.1cm}

The computations were done on a Sun Fire V890 Ultra Sparc IV
processor with 1200 MHz. Computation times are given in seconds,
rounded up. See Figure \ref{Figure: computation times}.
The running times give a clear ranking of the implementations: from
left to write the speed increases and in \emph{all} problems, the
presented Project-and-Lift algorithm wins significantly.

\begin{figure}[tbh]
{\footnotesize
\begin{center}
\[
\begin{array}{|l|r|r||r|r|r|r|r|r|}
  \hline
  \text{Problem} & \text{Vars.} & \text{GB size} & \text{Sing-blr}
  & \text{Sing-hs} & \text{CoCoA-t} & \text{CoCoA-ts} & \text{P\&L} &
  \FourTiTwo{\text{-gra}}\\
\hline
\hline
  333           & 27 &    110 & 30 &      4 &              1 &     1 &     1 & \\
  334           & 36 &    626 &  - &    197 &          3,024 &     5 &     1 & \\
  335           & 45 &  3,260 &  - & 23,700 &              - &   233 &    27 & \\
  344           & 48 &  7,357 &  - &      - &              - & 2,388 &   252 & \\
\hline
  \text{K}4     & 16 &     61 & 1 & 1 &         1 &     1 &     1 & \\
  \text{K}5     & 32 & 13,181 & - & - &    13,366 & 2,814 &   715 & \\
\hline
  \text{s-magic333}    & 66 &  1,664 & - & - & 35 &     55 &     3 & \\
  \text{grin}          &  8 &    214 & 4 & 4 &  1 &      1 &     1 & \\
\hline
\text{hppi}10 & 20 &  1,830 &  1,064 &    483 &        16 &     14 &     3 &      2 \\
\text{hppi}11 & 22 &  3,916 & 15,429 &  3,588 &       129 &     82 &    13 &     11 \\
\text{hppi}12 & 24 &  8,569 &      - & 43,567 &     1,534 &    554 &    60 &     51 \\
\text{hppi}13 & 26 & 16,968 &      - &      - &     8,973 &  4,078 &   290 &    259 \\
\text{hppi}14 & 28 & 34,355 &      - &      - &         - & 30,973 & 1,219 &  1,126\\
\hline
  \text{cuww}1  &  5 &      5 & - & - & - & - &     1 & \\
  \text{cuww}2  &  6 &     15 & - & - & - & - &     1 & \\
  \text{cuww}3  &  6 &     16 & - & - & - & - &     2 & \\
  \text{cuww}4  &  7 &      7 & - & - & - & - &     1 & \\
  \text{cuww}5  &  8 &     27 & - & - & - & - &     2 & \\
\hline
\end{array}
\]
\end{center}
}
\caption{Comparison of computing times.}
\label{Figure: computation times}
\end{figure}

The advantage of our Project-and-Lift algorithm is that it performs computations
in projected subspaces of $\Lattice$.  Thus, we obtain
comparably small intermediate sets during the computation.
Only the final iteration that deals with all variables reaches the true output
size. In contrast to this, the Saturation algorithm usually comes close to the
true output size already after the first saturation and then continues
computing with as many vectors.
See Figure \ref{Figure: Algorithm comparison}, for a comparison of intermediate
set sizes in each iteration for computing $3\times 4\times 4$ tables.
\input{Figure_Comparison}

Moreover, the Project-and-Lift algorithm, performs Gr\"obner basis computations
using a generating set, and thus can take full advantage of Criterion 2 which,
as computational experience shows, is extremely effective. In fact, we only
applied Criterion $2$ and $1$ (applied in that order) for the Project-and-Lift
algorithm since Criterion $3$ only slowed down the algorithm. However, for the
Saturation algorithm where we cannot apply Criterion $2$ fully, Criterion $3$
was very effective.  In this case, we applied Criterion $1$, then $3$, and then
$2$, in that order.

Note that in the knapsack problems \texttt{cuww1-cuww5} the initial set
$\sigma$ chosen by the Project-and-Lift Algorithm
\ref{Project-and-Lift algorithm} is empty. Thus, Algorithm
\ref{Project-and-Lift algorithm} simplifies to the Saturation
Algorithm \ref{Saturation algorithm}. In fact, only a single
saturation is necessary for each problem.

To us, the following observations were surprising.
\begin{itemize}
\item While Singular did not accept the inhomogeneous problems
\texttt{cuww1-cuww5} as input, CoCoA either could not solve them
or produced incorrect answers.
\item It is not clear why the CoCoA function Toric works well on
problems \texttt{hppi10-hppi14}, but runs badly on the table problems 334,
335, 344.
\item Problems \texttt{hppi10-hppi14} are in fact Graver basis
computations (see for example Chapter $14$ in \citet{Sturmfels:96}),
for which \FourTiTwo\ has the state-of-the-art algorithm and
implementation. Initially, it was a surprise to us that our
Project-and-Lift Algorithm \ref{Project-and-Lift algorithm} comes so
close to the speed of the state-of-the-art algorithm that computes
Graver bases directly \citep{Hemmecke:SymmGraver}. However, it turns
out that our Project-and-Lift Algorithm \ref{Project-and-Lift
algorithm} is an extension of the Project-and-Lift algorithm
presented in \citet{Hemmecke:SymmGraver} to lattice ideal
computations.
\end{itemize}

\bibliographystyle{elsart-harv}
\bibliography{bibliography}

\end{document}

%% file: Figure_Connected_Graph.tex
\begin{figure}[ht!]
\begin{center}
\framebox{ \setlength{\unitlength}{0.4pt}
\begin{picture}(310,260)(-20,-30)
    \put(260,210){(a)}
    \thicklines
    \put(0,0){\vector(0,1){220}}
    \multiputlist(0,-20)(50,0){0,1,2,3,4,5}
    \put(265,-20){$x_1$}
    \put(0,0){\vector(1,0){270}}
    \multiputlist(-20,0)(0,50){0,1,2,3,4}
    \put(-30,220){$x_2$}
    \matrixput(0,0)(50,0){6}(0,50){5}{\circle{7}}
    \dottedline{5}(0,50)(150,200)
    \dottedline{5}(0,100)(150,0)
    \dottedline{5}(100,0)(200,200)
    \dottedline{5}(150,200)(250,0)
    \put(100,100){\circle*{7}}
    \put(100,50){\circle*{7}}
    \put(50,100){\circle*{7}}
    \put(150,100){\circle*{7}}
    \put(100,150){\circle*{7}}
    \put(150,150){\circle*{7}}
    \put(150,200){\circle*{7}}
\end{picture}} \\
\vspace{0.3cm}
\framebox{ \setlength{\unitlength}{0.4pt}
\begin{picture}(310,260)(-20,-30)
    \put(260,210){(b)}
    \thicklines
    \put(0,0){\vector(0,1){220}}
    \multiputlist(0,-20)(50,0){0,1,2,3,4,5}
    \put(265,-20){$x_1$}
    \put(0,0){\vector(1,0){270}}
    \multiputlist(-20,0)(0,50){0,1,2,3,4}
    \put(-30,220){$x_2$}
    \matrixput(0,0)(50,0){6}(0,50){5}{\circle{7}}
    \dottedline{5}(0,50)(150,200)
    \dottedline{5}(0,100)(150,0)
    \dottedline{5}(100,0)(200,200)
    \dottedline{5}(150,200)(250,0)
    \put(100,100){\circle*{7}}
    \put(100,50){\circle*{7}}
    \put(50,100){\circle*{7}}
    \put(150,100){\circle*{7}}
    \put(100,150){\circle*{7}}
    \put(150,150){\circle*{7}}
    \put(150,200){\circle*{7}}


    \put(100,150){\line(1,0){50}}
    \put(100,100){\line(1,0){50}}
    \put(50,100){\line(1,0){50}}

    \put(150,100){\line(-1,1){50}}
    \put(100,50){\line(-1,1){50}}
\end{picture}}
\hspace{0.05cm}
\framebox{ \setlength{\unitlength}{0.4pt}
\begin{picture}(310,260)(-20,-30)
    \put(260,210){(c)}
    \thicklines
    \put(0,0){\vector(0,1){220}}
    \multiputlist(0,-20)(50,0){0,1,2,3,4,5}
    \put(265,-20){$x_1$}
    \put(0,0){\vector(1,0){270}}
    \multiputlist(-20,0)(0,50){0,1,2,3,4}
    \put(-30,220){$x_2$}
    \matrixput(0,0)(50,0){6}(0,50){5}{\circle{7}}
    \dottedline{5}(0,50)(150,200)
    \dottedline{5}(0,100)(150,0)
    \dottedline{5}(100,0)(200,200)
    \dottedline{5}(150,200)(250,0)
    \put(100,100){\circle*{7}}
    \put(100,50){\circle*{7}}
    \put(50,100){\circle*{7}}
    \put(150,100){\circle*{7}}
    \put(100,150){\circle*{7}}
    \put(150,150){\circle*{7}}
    \put(150,200){\circle*{7}}

    \put(100,150){\line(1,1){50}}
    \put(100,100){\line(1,1){50}}
    \put(50,100){\line(1,1){50}}
    \put(100,50){\line(1,1){50}}

    \put(100,150){\line(1,0){50}}
    \put(100,100){\line(1,0){50}}
    \put(50,100){\line(1,0){50}}

    \put(150,100){\line(-1,1){50}}
    \put(100,50){\line(-1,1){50}}
\end{picture}}
\end{center}
\caption{The set $\Feasible_{\Lattice,b}$ and the graphs $\Graph(\Feasible_{\Lattice,b},S)$ and
$\Graph(\Feasible_{\Lattice,b},S')$ projected onto the $(x_1,x_2)$-plane.}

\label{Figure: Connected Graph}
\end{figure}

%% file: Figure_Reduction_Path.tex
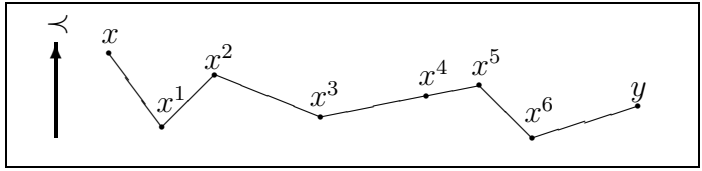
\begin{figure}[ht]
\begin{center}
\framebox{ \setlength{\unitlength}{0.4pt}
\begin{picture}(600,120)(-10,-10)

    \put(50,80){\circle*{5}}
    \put(100,10){\circle*{5}}
    \put(150,60){\circle*{5}}
    \put(250,20){\circle*{5}}
    \put(350,40){\circle*{5}}
    \put(400,50){\circle*{5}}
    \put(450,0){\circle*{5}}
    \put(550,30){\circle*{5}}

    \drawline(550,30)(450,0)(400,50)(350,40)(250,20)(150,60)(100,10)
        (50,80)

    \put(43,90){$x$}
    \put(543,40){$y$}

    \put(95,25){$x^1$}
    \put(140,65){$x^2$}
    \put(240,30){$x^3$}
    \put(343,50){$x^4$}
    \put(393,60){$x^5$}
    \put(443,15){$x^6$}

    \thicklines
    \put(0,0){\vector(0,1){90}}
    \put(-10,100){$\prec$}
\end{picture}}
\end{center}
\caption{Reduction path between $x$ and $y$.}

\label{Figure: Reduction Path}
\end{figure}

%% file: Figure_uv_Path.tex
\begin{figure}[ht!]
\begin{center}
\framebox{ \setlength{\unitlength}{0.2pt}
\begin{picture}(320,210)(-20,-20)

    \put(50,40){\circle*{10}}
    \put(150,140){\circle*{10}}
    \put(270,20){\circle*{10}}

    \drawline(50,40)(150,140)(270,20)

    \put(45,5){$x$}
    \put(135,155){$z$}
    \put(265,-15){$y$}

    \put(70,100){$u$}
    \put(220,85){$v$}

    \thicklines

    \put(0,0){\vector(0,1){160}}
    \put(-10,165){$\prec$}
\end{picture}}
\end{center}
\caption{A critical path for $(u,v)$ between $x$, $z$, and $y$.}

\label{Figure: uv Path}
\end{figure}
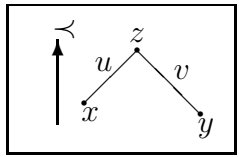

%% file: Figure_Replacing_Reduction_Path.tex
\begin{figure}[ht]
\begin{center}
\framebox{ \setlength{\unitlength}{0.35pt}
\begin{picture}(550,250)(0,0)

    \put(50,180){\circle*{5}}
    \put(100,110){\circle*{5}}
    \put(150,170){\circle*{5}}
    \put(250,210){\circle*{5}}
    \put(350,140){\circle*{5}}
    \put(400,50){\circle*{5}}
    \put(450,100){\circle*{5}}
    \put(500,160){\circle*{5}}

    \drawline(500,160)(450,100)(400,50)(350,140)(250,210)(150,170)(100,110)
        (50,180)

    \put(43,190){$x^1$}
    \put(493,170){$x^8$}

    \put(95,125){$x^2$}
    \put(137,175){$x^3$}
    \put(235,215){$x^4$}
    \put(343,150){$x^5$}
    \put(393,65){$x^6$}
    \put(443,115){$x^7$}

    \put(300,180){$v$}
    \put(190,195){$u$}

    \put(200,100){\circle*{5}}
    \put(250,60){\circle*{5}}
    \put(300,80){\circle*{5}}
    \dottedline{5}(350,140)(300,80)(250,60)(200,100)(150,170)
    \put(135,130){$\overline{x}^0$}
    \put(180,65){$\overline{x}^1$}
    \put(235,30){$\overline{x}^2$}
    \put(295,50){$\overline{x}^3$}
    \put(330,95){$\overline{x}^4$}

    \thicklines
    \put(0,0){\vector(0,1){220}}
    \put(-10,235){$\prec$}
\end{picture}}
\end{center}
\caption{Replacing a critical path by a reduction path}

\label{Figure: Replacing Reduction Path}
\end{figure}
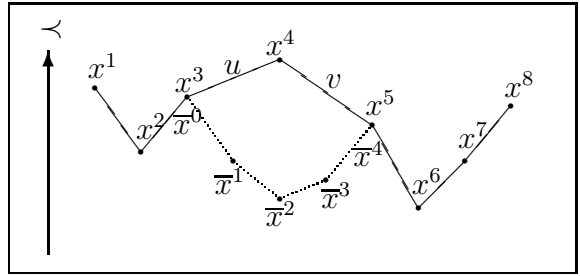

%% file: Figure_Criterion_1.tex
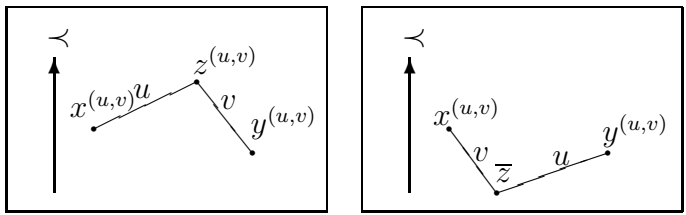
\begin{figure}[ht!]
\begin{center}
\framebox{ \setlength{\unitlength}{0.3pt}
\begin{picture}(330,210)(-10,0)

    \put(50,80){\circle*{5}}
    \put(180,140){\circle*{5}}
    \put(250,50){\circle*{5}}

    \drawline(50,80)(180,140)(250,50)

    \put(20,95){$x^{(u,v)}$}
    \put(175,150){$z^{(u,v)}$}
    \put(248,75){$y^{(u,v)}$}

    \put(100,120){$u$}
    \put(210,105){$v$}

    \thicklines

    \put(0,0){\vector(0,1){170}}
    \put(-10,185){$\prec$}
\end{picture}}
\hspace{5pt} \framebox{ \setlength{\unitlength}{0.3pt}
\begin{picture}(330,210)(-10,0)

    \put(50,80){\circle*{5}}
    \put(110,0){\circle*{5}}
    \put(250,50){\circle*{5}}

    \drawline(50,80)(110,0)(250,50)

    \put(30,85){$x^{(u,v)}$}
    \put(108,13){$\overline{z}$}
    \put(243,65){$y^{(u,v)}$}

    \put(180,35){$u$}
    \put(80,40){$v$}

    \thicklines

    \put(0,0){\vector(0,1){170}}
    \put(-10,185){$\prec$}
\end{picture}}
\end{center}
\caption{Criterion $1$.}

\label{Figure: Criterion 1}

\end{figure}

%% file: Figure_Criterion_2.tex
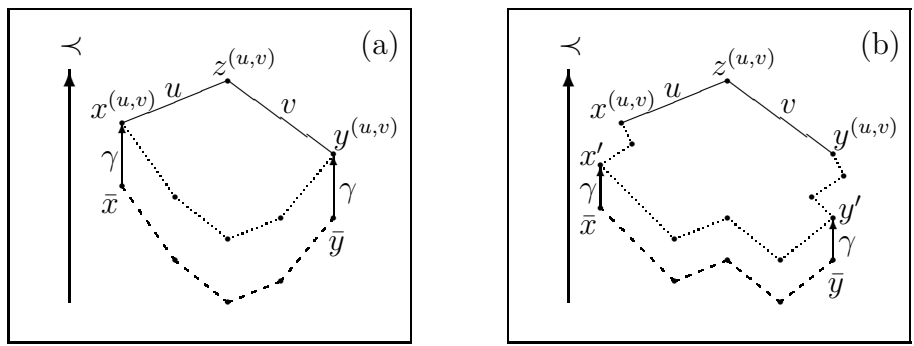
\begin{figure}[ht!]
\begin{center}
\framebox{ \setlength{\unitlength}{0.4pt}
\begin{picture}(325,280)(-20,-20)
    \put(275, 235){(a)}

    \put(50,170){\circle*{5}}
    \put(150,210){\circle*{5}}
    \put(250,140){\circle*{5}}

    \drawline(250,140)(150,210)(50,170)

    \put(20,175){$x^{(u,v)}$}
    \put(135,215){$z^{(u,v)}$}
    \put(250,150){$y^{(u,v)}$}

    \put(200,180){$v$}
    \put(90,195){$u$}

    \put(100,100){\circle*{5}}
    \put(150,60){\circle*{5}}
    \put(200,80){\circle*{5}}
    \dottedline{5}(250,140)(200,80)(150,60)(100,100)(50,170)

    \put(50,110){\circle*{5}}
    \put(100,40){\circle*{5}}
    \put(150,0){\circle*{5}}
    \put(200,20){\circle*{5}}
    \put(250,80){\circle*{5}}
    \dashline{5}(250,80)(200,20)(150,0)(100,40)(50,110)

    \put(30,85){$\bar{x}$}
    \put(245,50){$\bar{y}$}

    \put(250,80){\vector(0,1){60}}
    \put(50,110){\vector(0,1){60}}

    \put(30,130){$\gamma$}
    \put(255,100){$\gamma$}

    \thicklines

    \put(0,0){\vector(0,1){220}}
    \put(-10,235){$\prec$}
\end{picture}}
\hspace{1cm}
\framebox{ \setlength{\unitlength}{0.4pt}
\begin{picture}(325,280)(-20,-20)
    \put(275, 235){(b)}

    \put(50,170){\circle*{5}}
    \put(150,210){\circle*{5}}
    \put(250,140){\circle*{5}}

    \drawline(250,140)(150,210)(50,170)

    \put(20,175){$x^{(u,v)}$}
    \put(135,215){$z^{(u,v)}$}
    \put(250,150){$y^{(u,v)}$}

    \put(200,180){$v$}
    \put(90,195){$u$}

    \put(60,150){\circle*{5}}
    \dottedline{5}(50,170)(60,150)(30,130)

    \put(230,100){\circle*{5}}
    \put(260,120){\circle*{5}}
    \dottedline{5}(250,80)(230,100)(260,120)(250,140)

    \put(10,130){$x'$}
    \put(255,80){$y'$}

    \put(30,130){\circle*{5}}
    \put(100,60){\circle*{5}}
    \put(150,80){\circle*{5}}
    \put(200,40){\circle*{5}}
    \put(250,80){\circle*{5}}
    \dottedline{5}(250,80)(200,40)(150,80)(100,60)(30,130)

    \put(30,90){\circle*{5}}
    \put(100,20){\circle*{5}}
    \put(150,40){\circle*{5}}
    \put(200,0){\circle*{5}}
    \put(250,40){\circle*{5}}
    \dashline{5}(250,40)(200,0)(150,40)(100,20)(30,90)

    \put(10,65){$\bar{x}$}
    \put(245,10){$\bar{y}$}

    \put(30,90){\vector(0,1){40}}
    \put(250,40){\vector(0,1){40}}

    \put(10,100){$\gamma$}
    \put(255,50){$\gamma$}

    \thicklines

    \put(0,0){\vector(0,1){220}}
    \put(-10,235){$\prec$}
\end{picture}}
\end{center}
\caption{Criterion $2$.}

\label{Figure: Criterion 2}
\end{figure}

%% file: Figure_Criterion_3.tex
\begin{figure}[ht!]
\begin{center}
\framebox{ \setlength{\unitlength}{0.35pt}
\begin{picture}(300,210)(-10,0)
    \put(250,185){(a)}

    \put(50,80){\circle*{5}}
    \put(170,160){\circle*{5}}
    \put(250,50){\circle*{5}}
    \put(130,35){\circle*{5}}
    
    \drawline(50,80)(170,160)(250,50)
    \drawline(130,35)(170,160)

    \put(20,50){$x^{(u,v)}$}
    \put(165,170){$z^{(u,v)}$}
    \put(225,25){$y^{(u,v)}$}
    \put(125,10){$\overline{z}$}

    \put(90,120){$u$}
    \put(220,105){$v$}
    \put(155,95){$w$}

    \thicklines

    \put(0,0){\vector(0,1){170}}
    \put(-10,185){$\succ$}
\end{picture}}
\hspace{1cm}
\framebox{ \setlength{\unitlength}{0.35pt}
\begin{picture}(300,210)(-10,0)
    \put(250,185){(b)}

    \put(50,80){\circle*{5}}
    \put(170,160){\circle*{5}}
    \put(250,50){\circle*{5}}
    \put(130,90){\circle*{5}}
    
    \drawline(50,80)(170,160)(250,50)
    \drawline(130,90)(170,160)

    \put(20,50){$x^{(u,v)}$}
    \put(165,170){$z^{(u,v)}$}
    \put(225,25){$y^{(u,v)}$}
    \put(125,65){$\overline{z}$}

    \put(90,120){$u$}
    \put(220,105){$v$}
    \put(155,115){$w$}

    \thicklines

    \put(0,0){\vector(0,1){170}}
    \put(-10,185){$\succ$}
\end{picture}}
\end{center}
\caption{Criterion $3$.}

\label{Figure: Criterion 3}
\end{figure}
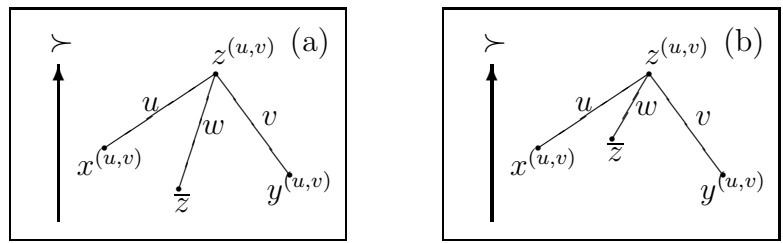

%% file: Figure_Comparison.tex
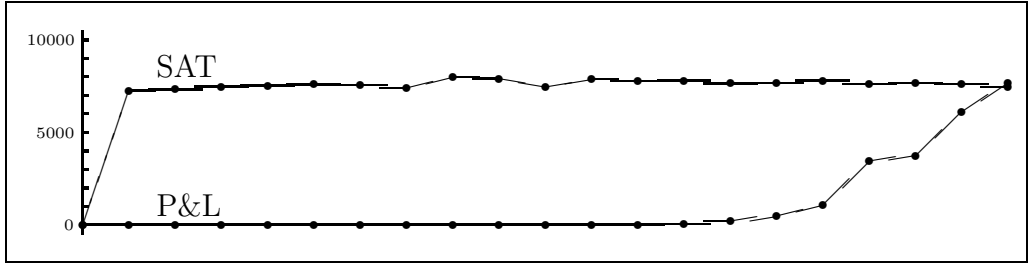
\begin{figure}[ht]
\begin{center}
\framebox{ \setlength{\unitlength}{0.7pt}
\begin{picture}(520,120)(-20,-10)
    \drawline   (0,0.18)
                (25,72.52)
                (50,73.23)
                (75,74.73)
                (100,75.12)
                (125,75.93)
                (150,75.73)
                (175,73.91)
                (200,79.92)
                (225,79.07)
                (250,74.78)
                (275,78.77)
                (300,77.76)
                (325,77.83)
                (350,76.54)
                (375,76.66)
                (400,77.83)
                (425,76.13)
                (450,76.71)
                (475,76.01)
                (500,74.48)

    \put(0,0.18){\circle*{4}}
    \put(25,72.52){\circle*{4}}
    \put(50,73.23){\circle*{4}}
    \put(75,74.73){\circle*{4}}
    \put(100,75.12){\circle*{4}}
    \put(125,75.93){\circle*{4}}
    \put(150,75.73){\circle*{4}}
    \put(175,73.91){\circle*{4}}
    \put(200,79.92){\circle*{4}}
    \put(225,79.07){\circle*{4}}
    \put(250,74.78){\circle*{4}}
    \put(275,78.77){\circle*{4}}
    \put(300,77.76){\circle*{4}}
    \put(325,77.83){\circle*{4}}
    \put(350,76.54){\circle*{4}}
    \put(375,76.66){\circle*{4}}
    \put(400,77.83){\circle*{4}}
    \put(425,76.13){\circle*{4}}
    \put(450,76.71){\circle*{4}}
    \put(475,76.01){\circle*{4}}
    \put(500,74.48){\circle*{4}}

    \drawline   (0,0.18)
                (25,0.18)
                (50,0.18)
                (75,0.18)
                (100,0.18)
                (125,0.18)
                (150,0.18)
                (175,0.18)
                (200,0.18)
                (225,0.18)
                (250,0.18)
                (275,0.18)
                (300,0.18)
                (325,0.39)
                (350,2.06)
                (375,4.57)
                (400,10.57)
                (425,34.31)
                (450,37.36)
                (475,61.06)
                (500,76.78)

    \put(0,0.18){\circle*{4}}
    \put(25,0.18){\circle*{4}}
    \put(50,0.18){\circle*{4}}
    \put(75,0.18){\circle*{4}}
    \put(100,0.18){\circle*{4}}
    \put(125,0.18){\circle*{4}}
    \put(150,0.18){\circle*{4}}
    \put(175,0.18){\circle*{4}}
    \put(200,0.18){\circle*{4}}
    \put(225,0.18){\circle*{4}}
    \put(250,0.18){\circle*{4}}
    \put(275,0.18){\circle*{4}}
    \put(300,0.18){\circle*{4}}
    \put(325,0.39){\circle*{4}}
    \put(350,2.06){\circle*{4}}
    \put(375,4.57){\circle*{4}}
    \put(400,10.57){\circle*{4}}
    \put(425,34.31){\circle*{4}}
    \put(450,37.36){\circle*{4}}
    \put(475,61.06){\circle*{4}}
    \put(500,76.78){\circle*{4}}

    \thicklines
    \put(0,-5){\line(0,1){110}}
    \put(-10,-3){\tiny{0}}
    \put(-25,47){\tiny{5000}}
    \put(-30,97){\tiny{10000}}
    \multiput(0,0)(0,10){11}{\line(1,0){3}}
    \put(40,5){P\&L}
    \put(40,80){SAT}
\end{picture}}
\end{center}
\caption{Comparison of intermediate set sizes in each iteration.}

\label{Figure: Algorithm comparison}
\end{figure}